\documentclass[11pt]{article} 
\input epsf
\input amssym.def
\input amssym.tex

\usepackage{graphicx}
\usepackage{caption2}
\usepackage{thm}

\newcommand{\beq}{\begin{equation}}
\newcommand{\CAC}{{\cal C}}
\newcommand{\CAH}{{\cal H}}

\newcommand{\diam}{\mbox{diam}}
\newcommand{\eeq}{\end{equation}}
\newcommand{\gnaive}{\widetilde{1/\gamma}_{\eta}}
\newcommand{\gols}{\widehat{1/\gamma}_{\eta}}
\newcommand{\gr}{{\rm\mbox{gr}}}
\newcommand{\hnaive}{\widetilde{1/\upsilon}_{\eta}}
\newcommand{\hols}{\widehat{1/\upsilon}_{\eta}}
\newcommand{\naive}{\widetilde{\upsilon}_{\eta}}
\newcommand{\ols}{\widehat{\upsilon}_{\eta}}
\newcommand{\osc}{{\rm\mbox{Osc}}}
\newcommand{\prodi}{\mathop{{\lower3pt\hbox{\epsfxsize=7pt\epsfbox{pi.ps}}}}}

\newcommand{\RR}{\Bbb R}
\newcommand{\sas}{S\alpha S}

\begin{document}

\title{Estimating the $p$-variation index of a sample function:\\
An application to financial data set}
%\vskip 1pc
\author{Rimas Norvai\v sa\thanks{This research was partially supported 
by NSERC Canada grant number 203232-98 at York University and by
Lithuanian State Science and Studies Foundation Grant K014.}\\
Vilnius
\and Donna Mary Salopek\thanks{This research was partially supported 
by NSERC Canada grant number 203232-98 at York University}\\
Toronto}
\date{\today}

%\maketitle

\noindent
{\bf\LARGE Estimating the $p$-variation index of a sample function:\\
An application to financial data set}

\vskip 2pc
\noindent
RIMAS NORVAI\v SA\footnote{This research was partially supported 
by NSERC Canada grant number 203232-98 at York University and by
Lithuanian State Science and Studies Foundation Grant K014.}
\hfill norvaisa@ktl.mii.lt

\noindent
{\em\small Institute of Mathematics and Informatics, Akademijos 4, 
Vilnius 2600, Lithuania}

\vskip 1pc
\noindent
DONNA MARY SALOPEK\footnote{This research was partially supported 
by NSERC Canada grant number 203232-98 at York University.}
\hfill dsalopek@mathstat.yorku.ca

\noindent
{\em\small York University, Department of Mathematics and Statistics,
4700 Keele Street,} 

\noindent
{\em\small Toronto, Ontario M3J 1P3 Canada}

\vskip 1pc
\noindent
{\small {\bf Abstract.}
In this paper we apply a real analysis approach to test continuous
time stochastic models of financial mathematics.
Specifically, fractal dimension estimation methods are applied to
statistical analysis of continuous time stochastic processes.
To estimate a roughness of a sample function we modify a box-counting 
method typically used in estimating fractal dimension of a graph
of a function. 
Here the roughness of a function $f$ is defined as the infimum of numbers
$p>0$ such that $f$ has bounded $p$-variation, which we call the 
$p$-variation index of $f$.
The method is also tested on estimating the exponent $\alpha\in [1,2]$ of a 
simulated symmetric $\alpha$-stable process, and on estimating the Hurst 
exponent $H\in (0,1)$ of a simulated fractional Brownian motion.}



\section{Introduction: the $p$-variation index and financial modeling}

In financial mathematics, a simplest continuous time  
model assumes that a stock price, or other financial asset,
is a stochastic process $P=\{P(t)\colon\,0\leq t\leq T\}$ satisfying
the relation 
\beq\label{BSM-model}
P(t)=1+(I)\int_0^tP\,dB,\qquad 0\leq t\leq T,
\eeq
which is the It\^o integral equation with respect to a standard Brownian 
motion $B=\{B(t)\colon\,t\geq 0\}$.
Equation (\ref{BSM-model}) is usually written in the form of
a stochastic differential equation, or simply by giving its solution
$P(t)=\exp\{B(t)-t/2\}$, $0\leq t\leq T$.
In the financial literature this is known as the Black-Scholes-Merton 
stock price model.
Its pertinence is backed-up by the assumption that increments of the log 
transform of a stock price are independent and normally distributed, which
is known as the strong form of Random Walk Hypothesis.
In agreement with a relaxed form or an alternative form of the Random Walk 
Hypothesis, the Brownian motion $B$ in equation (\ref{BSM-model}) can be 
replaced by a more general stochastic process $X$, 
and the linear It\^o integral equation
(\ref{BSM-model}) by a different integral equation with respect to $X$.
In this paper, the stochastic process $X$ is called the {\em return process},
and the unique solution of an integral equation with respect to $X$ 
is called the {\em stock price process} $P$. 

The mainstream econometric analysis of continuous time financial models is
to test different hypotheses about an integral equation describing a stock 
price process $P$, or to test various parameters of a distribution of a return 
process $X$ (see e.g.\ Section 9.3 in Campbell, Lo and MacKinlay, 1997).
In this paper we attempt to test the degree of roughness of a
return process $X$.
The legitimacy of such an endeavor is based on the fact that
the support of a distribution of a stochastic processes is
a particular class of functions.
That is,  a suitable class of functions contains almost all
sample functions  of a stochastic process.
More specifically, if $X$ is a regular enough stochastic process
defined on a given probability space $(\Omega,{\cal F},\Pr)$, then for a 
suitable class of functions $F$, a {\em sample function} 
$X(\omega)=\{X(t,\omega)\colon\, 0\leq t\leq T\}$ belongs to $F$ 
for almost all $\omega\in\Omega$.
Often $F$ can be taken as a proper subspace of the space of all 
continuous functions on $[0,T]$, or a proper subspace of the (Skorohod) 
space of all regulated and right continuous functions on $[0,T]$.
For example, the support of the distribution of a standard Brownian
motion is a subset of a class of functions having the order
of H\"older continuity strictly bigger than $1/2$.
However, the H\"older continuity is not applicable to characterize 
sample functions of a L\'evy stochastic process without a Gaussian component,
because almost every sample function of such a process is discontinuous. 
A simple example of a L\'evy process is a symmetric $\alpha$-stable
stochastic process, or a $\sas$ stochastic process, $X_{\alpha}$ with 
the exponent $\alpha\in (0,2]$.
The case $\alpha =2$ gives the only Gaussian component, a standard Brownian 
motion, that is, $X_2=B$.
In general, if $X_{\alpha}$ is a $\sas$ stochastic process with 
$\alpha\in (0,2]$,  then almost all sample functions of $X_{\alpha}$ 
have bounded $p$-variation for each $p>\alpha$, and have infinite
$p$-variation if $0<p\leq\alpha$.
This fact was known since the time P.\ L\'evy first introduced 
this process around the beginning of thirties.

To recall the property of boundedness of $p$-variation, let
$f$ be a real-valued function on an interval $[0,T]$.
For a number $0<p<\infty$, let 
$$
s_p(f;\kappa):=\sum_{i=1}^n|f(t_i)-f(t_{i-1})|^p, 
$$
where $\kappa=\{t_i\colon\,i=0,\dots,n\}$ is a partition of $[0,T]$,
that is $0=t_0<t_1<\cdots <t_n=T$.
The {\em $p$-variation} of $f$ is defined by
$$
v_p(f):=v_p(f;[0,T]):=\sup\big\{s_p(f;\kappa)\colon\,
\kappa\,\,\mbox{is a partition of}\,\, [0,T]\big\}.
$$ 
The function $f$ has bounded $p$-variation on $[0,T]$ if $v_p(f)<\infty$.
By H\"older's inequality it follows that if $v_p(f)<\infty$ and
$q>p$, then $v_q(f)<\infty$.
The number
$$
\upsilon (f;[0,T]):=\inf\big\{p>0\colon\,v_p(f;[0,T])<\infty\big\}
=\sup\big\{p>0\colon\,v_p(f;[0,T])=\infty\big\}
$$ 
is called the {\em $p$-variation index} of $f$.
For a regular enough stochastic process $X$, $\upsilon (X)(\omega):=
\upsilon (X(\cdot,\omega))$, $\omega\in\Omega$, is a random variable
which we call the $p$-variation index of $X$.
In fact, for all stochastic processes discussed in this paper, their
$p$-variation indices are known to be constants and we seek
to estimate these constants.   

Now for a $\sas$ stochastic process
$X_{\alpha}$, we can restate its sample regularity by saying that 
its $p$-variation index $\upsilon (X_{\alpha})=\alpha$ almost surely.
This is a special case of the following more general fact.

\begin{exmp}\label{exmp1}
Let $X$ be a homogeneous L\'evy stochastic process with the L\'evy
measure $\nu$, which is a $\sigma$-finite Borel measure on $\RR\setminus \{0\}$
such that
$$
\int_{\RR\setminus\{0\}}\min\{1, |x|^2\}\,\nu (dx)<+\infty.
$$
The {\em Blumenthal-Getoor index} $\beta_X$ of $X$ is defined by
$$
\beta_X :=\inf\big\{ \alpha>0\colon\,\int_{\RR\setminus\{0\}}\min\{1,
|x|^{\alpha}\}\,\nu (dx)<+\infty\big\}.
$$
Note that $0<\beta_X\leq 2$.
If $X$ has no Gaussian part, then for any $0<T<\infty$,
$$ 
\upsilon(X;[0,T])=\beta_X,\qquad\mbox{almost surely.} 
$$
This follows from Theorems 4.1 and 4.2 of Blumenthal and Getoor (1961),
and from Theorem 2 of Monroe (1972).
\end{exmp}

A stock price model having as a return a L\'evy process without
a Gaussian part is a common alternative to the Black-Scholes-Merton model.
Another popular alternative is a fractional Brownian
motion $B_H$ with the Hurst exponent $H\in (0,1)$, where
$B_H$ with $H=1/2$ is a standard Brownian motion. 
The $p$-variation index of a fractional Brownian motion
$\upsilon (B_H;[0,T])=1/H$ almost surely.
This is the special case of the following fact.

\begin{exmp}\label{exmp2}
Let $X=\{X(t)\colon\,t\geq 0\}$ be a Gaussian stochastic process
with stationary increments and continuous in quadratic mean.
Let $\sigma_X$ be the incremental variance of $X$ defined by
$\sigma_X(u)^2:=E[(X(t+u)-X(t))^2]$ for $t, u\geq 0$.
Let
$$
\gamma_{\ast}:=\inf\big\{\gamma >0\colon\,\lim_{u\downarrow 0}
\frac{u^{\gamma}}{\sigma_X(u)}=0\big\}\quad\mbox{and}\quad
\gamma^{\ast}:=\sup\big\{\gamma >0\colon\,\lim_{u\downarrow 0}
\frac{u^{\gamma}}{\sigma_X(u)}=+\infty\big\}.
$$
Then $0\leq\gamma_{\ast}\leq\gamma^{\ast}\leq +\infty$.
If $\gamma_{\ast}=\gamma^{\ast}$, then we say that $X$ has
an {\em Orey index} $\gamma_X:=\gamma_{\ast}=\gamma^{\ast}$.
Furthermore, if $X$ has an Orey index $\gamma_X\in (0,1)$, then for any
$0<T<+\infty$,
\beq\label{var-of-Orey}
\upsilon (X;[0,T])=1/\gamma_X,\qquad\mbox{almost surely.}
\eeq
This follows from the fact that almost all sample functions of
$X$ obey a uniform H\"older condition with exponent $\gamma
<\gamma_X$ (see Section 9.4 of Cramer and Leadbetter, 1967)
and from the inequality of Berman (1969) connecting the
$p$-variation with the Fourier transform of local times of $X$.
Relation (\ref{var-of-Orey}) also follows from the characterization
of the $p$-variation index for arbitrarily Gaussian processes due to
Jain and Monrad (1983).
\end{exmp}

From the point of view of a statistical time series analysis,
estimation of the $p$-variation index in the above two examples offer 
a new perspective to analyzing financial data sets.
For instance, a symmetric $\alpha$-stable process and a fractional 
Brownian motion with the Hurst exponent $H$, both have the same $p$-variation
index in the case $\alpha=1/H\in (1,2)$. 
However the latter has exponentially small tails, while the former 
has not even the second moment.
These two examples are extensions of the Black-Scholes-Merton model
(that is when $\alpha=1/H=2$) into two different directions.
The Orey index, and so the $p$-variation index by relation (\ref{var-of-Orey}),
have already been estimated in the paper Norvai\v sa and Salopek (2000).
They used two estimators based on the result of Gladyshev (1961).  
The estimators of the present paper can be applied under much less
restrictive hypotheses about stock price returns, and helps to reconcile
the two divergent directions of theoretical analyses of financial markets. 

\section{The oscillation $\eta$-summing index and related estimators}

In this section we describe a method of estimating the $p$-variation
index of a function based on existence of the metric entropy index
(or the box-counting dimension) of its graph.
Let $f$ be a real-valued function defined on an interval $[0,T]$.
Let $\eta=\{N_m\colon\,m\geq 1\}$ be a sequence of strictly increasing
positive integers.
With $\eta$ one can associate a sequence $\{\lambda (m)\colon\,m\geq 1\}$
of partitions $\lambda (m)=\{iT/N_m\colon\,i=0,\dots,N_m\}$ of $[0,T]$ into 
subintervals $\Delta_{i,m}:=[(i-1)T/N_m,iT/N_m]$, $i=1,\dots,N_m$, 
all having the same length $T/N_m$.
For each $m\geq 1$, let
\beq\label{Q(f)}
Q(f;\lambda (m))=\sum_{i=1}^{N_m}\osc (f;\Delta_{i,m}),
\eeq
where for a subset $A\subset [0,T]$,
$$
\osc (f;A):=\sup\big\{|f(t)-f(s)|\colon\,s, t\in A\big\}
=\sup_{t\in A}f(t)-\inf_{s\in A}f(s).
$$
The sequence $Q_{\eta}(f):=\{Q(f;\lambda (m))\colon\,m\geq 1\}$ will be called
the {\em oscillation $\eta$-summing} sequence.
For a bounded non-constant function $f$ on $[0,T]$, and any sequence $\eta$ 
as above, let
$$
\delta_{\eta}^{-}(f):=
\liminf_{m\to\infty}\frac{\log Q(f;\lambda (m))/N_m}{\log (1/N_m)}
\quad\mbox{and}\quad\delta_{\eta}^{+}(f)
:=\limsup_{m\to\infty}\frac{\log Q(f;\lambda (m))/N_m}{\log (1/N_m)}.
$$
Then we have
\beq\label{Qin[1,2]}
0\leq\delta_{\eta}^{-}(f)\leq\delta_{\eta}^{+}(f)\leq 1.
\eeq
The lower bound follows from the bound $Q(f;\lambda (m))/N_m\leq 
\osc (f;[0,T])<\infty$, which is valid for each $m\geq 1$.
The upper bound holds because $Q(f;\lambda (m))\geq C/2>0$ for 
all sufficiently large $m\geq 1$.
Indeed, if $f$ is continuous, then $C=v_1(f;[0,T])$. 
Otherwise $f$ has a jump at some $t\in [0,T]$, so that $C$ can be taken to be
a saltus at $t$ if $t\not\in\lambda (m)$ for all sufficiently large $m$,
or $C$ can be taken to be a one-sided non-zero saltus at $t$ if 
$t\in\lambda (m)$ for infinitely many $m$. 
Instead of relation (\ref{Qin[1,2]}), a sharp one-sided bound is given by
Lemma \ref{bc-pv} in Appendix A.

If $f$ has bounded variation, then for each $\eta$,
\beq\label{bv}
\delta_{\eta}^{-}(f)=\delta_{\eta}^{+}(f)=1.
\eeq
Indeed, since $Q(f;\lambda (m))\leq v_1(f;[0,T])<\infty$,
we have $\log Q(f;\lambda (m))/N_m\geq \log v_1(f;[0,T])/N_m$
for all sufficiently large $m$.
Thus $\delta_{\eta}^{-}(f)\geq 1$, and equalities in (\ref{bv}) follow from 
relation (\ref{Qin[1,2]}).

\begin{defn}\label{OS-index}
Let $f$ be a non-constant real-valued function on $[0,T]$,
and let $\eta=\{N_m\colon\,m\geq 1\}$ be a sequence of strictly increasing
positive integers.
If $\delta_{\eta}^{-}(f)=\delta_{\eta}^{+}(f)$, then we say that $f$ has
the {\em oscillation $\eta$-summing index} $\delta_{\eta}(f)$ 
and is defined by
\beq\label{1OS-index}
\delta_{\eta}(f):=\delta_{\eta}^{-}(f)=\delta_{\eta}^{+}(f)
=\lim_{m\to\infty}\frac{\log Q(f;\lambda (m))/N_m}{\log (1/N_m)}. 
\eeq
\end{defn}

Next we give a sufficient condition for exitence of the oscillation
$\eta$-summing index for any $\eta$.
Let $E$ be a nonempty bounded subset in a plane $\RR^2$, and let 
$N(E;\epsilon)$, $\epsilon >0$, be the minimum number of closed balls of 
diameter $\epsilon$  required to cover $E$.
The {\em lower} and {\em upper metric entropy indices} of the set $E$ are 
defined respectively by
$$
\Delta^{-}(E):=\liminf_{\epsilon\downarrow 0}\frac{\log N(E;\epsilon)}
{\log (1/\epsilon)}\quad\mbox{and}\quad
\Delta^{+}(E):=\limsup_{\epsilon\downarrow 0}\frac{\log N(E;\epsilon)}
{\log (1/\epsilon)}.
$$
If $\Delta^{-}(E)=\Delta^{+}(E)$, then the common value denoted by
$\Delta (E)$ is called the {\em metric entropy index} of the set $E$.
In the actual calculations of the metric entropy index, it is often simpler
to replace closed balls by squares (boxes) of a grid (cf.\ Lemma \ref{boxes}
below).
Therefore in fractal analysis, $\Delta (E)$ is also known as 
the box-counting dimension.

The proof of the following theorem is given in Appendix A.

\begin{thm}\label{main}
Let $f$ be a regulated non-constant function on $[0,T]$ with the $p$-variation
index $\upsilon (f)$.
If the metric entropy index of the graph $\gr (f)$  of $f$ is defined
and
\beq\label{2main}
\Delta (\gr (f))=2-1/(1\vee\upsilon (f)), 
\eeq
then for any sequence $\eta$, $f$ has the oscillation $\eta$-summing index 
\beq\label{index=var}
\delta_{\eta}(f)=1/(1\vee\upsilon (f)).
\eeq 
\end{thm}

Essentially, relation (\ref{2main}) is the  lower bound condition
on the metric entropy index because the following always holds. 

\begin{prop}\label{3main}
For a regulated function $f$ on $[0,T]$,
$\Delta^{+} (\gr (f))\leq 2-1/(1\vee\upsilon (f))$.
\end{prop}

The proof is similar to the proof of Theorem \ref{main}
and is also given in Appendix A.

The oscillation $\eta$-summing index is a slightly modified concept
of a real box-counting method introduced by Carter, Cawley and Mauldin (1988).
Independently, Dubuc, Quiniou, Roques-Carmes, Tricot and Zucker 
(1989) arrived at essentially the same notion, but called it
the variation method.
Both papers applied the new index to estimate 
the fractal dimension of several continuous functions whose dimension
is known, and found the new algorithm superior over several other
fractal dimension estimation methods (see also Section 6.2 in 
Cutler, 1993, for further discussion on this).

\paragraph*{Oscillation $\eta$-summing estimators.}
Let $f$ be a real-valued function defined on $[0,1]$, and let 
$\eta=\{N_m\colon\,m\geq 1\}$ be a sequence of strictly increasing
positive integers.
Let $\{u_1,\dots,u_N\}\subset [0,1]$ be a set of points such that
for some integer $M$,
\beq\label{1ose}
\cup_{m=1}^M\lambda (m)=\{u_1,\dots,u_N\},\qquad
\mbox{where $\lambda (m):=\{i/N_m\colon\,i=0,\dots,N_m\}$.}
\eeq
Given a finite set of values $\{f(u_1),\dots,f(u_N)\}$,
we want to estimate the $p$-variation
index $\upsilon (f)$.
To achieve this, for each $m\in\{1,\dots,M\}$, let
\beq\label{Q(m)}
Q(m):=\sum_{i=1}^{N_m}\Big [\max_{u_k\in\Delta_{i,m}}\big\{ f(u_k)\big\}
-\min_{u_k\in\Delta_{i,m}}\big\{ f(u_k)\big\}\Big ],
\eeq
where $\Delta_{i,m}=[(i-1)/N_m,i/N_m]$.
For large enough $M$, the finite set $\{Q(m)\colon\,m=1,\dots,M\}$
may be considered as an approximation to the oscillation $\eta$-summing 
sequence $Q_{\eta}(f)$ defined by (\ref{Q(f)}).
For $m=1,\dots,M$, let 
\beq\label{r_m}
r(m):=\frac{\log_2 1/N_m}{\log_2 Q(m)/N_m}
=\frac{\log_2 N_m}{\log_2 N_m/Q(m)}.
\eeq
Relation (\ref{index=var}) suggests that the set $\{r(m)\colon\,
m=1,\dots,M\}$ may be used to estimate the $p$-variation index
$\upsilon (f)$.

\begin{defn}\label{ose}
Let $\eta=\{N_m\colon\,m\geq 1\}$ be a sequence of strictly increasing 
positive integers, let $\{u_1,\dots,u_N\}\subset [0,1]$ be such that 
(\ref{1ose}) holds for some integer $M$, and let $\{f(u_1),\dots,f(u_N)\}$ 
be a set of known values of a real-valued function $f$ on $[0,1]$.
We will say that $\naive (f):=r(M)$ is the {\em naive oscillation
$\eta$-summing estimator} of $1\vee \upsilon (f)$ based on
$\eta_M:=\{N_m\colon\,1\leq m\leq M\}$.
Letting $x_m:=\log_2 (N_m/Q(m))$ and $\bar x:=M^{-1}\sum_{m=1}^Mx_m$,
we will say that
$$
\ols (f):=\frac{\sum_{m=1}^M(x_m-\bar x)\log_2 N_m}{\sum_{m=1}^M
(x_m-\bar x)^2}
$$ 
is the {\em OLS oscillation $\eta$-summing estimator} of $1\vee\upsilon (f)$
based on $\eta_M:=\{N_m\colon\,1\leq m\leq M\}$.
The estimators $\naive$ and $\ols$ will be called the OS estimators, and
the estimation either by $\naive$ or by $\ols$ will be called the OS 
estimation.
\end{defn}

Relation (\ref{index=var}) alone, if it holds for a function $f$ and 
a sequence $\eta$, does not imply that the two estimators will converge 
to $\upsilon (f)$ as $N\to\infty$, and so as $M\to\infty$ by relation
(\ref{1ose}).
If $\upsilon (f)<\infty$, then $f$ is a regulated function on $[0,1]$,
that is, there exist the limits $f(t+):=\lim_{u\downarrow t}f(u)$ for each 
$t\in [0,1)$ and $f(s-):=\lim_{u\uparrow s}f(u)$ for each $s\in (0,1]$.
Assuming that $f$ is regulated and either right- or left-continuous,
then $\osc (f;A)$ is the same as $\osc (f;A\cap U)$ for a countable
and dense subset $U\subset [0,1]$ and any subset $A\subset [0,1]$.
For such a function $f$, one can show that the naive estimator $\naive (f)$
will approach $\upsilon (f)$ as the set $\{u_1,\dots,u_N\}$ will 
increase to $\cup_m\lambda (m)$.
For sample functions of a stationary Gaussian stochastic process $X$,
Hall and Wood (1993) showed that the two estimators corresponding to
the reciprocal of relation (\ref{r_m}) converge to $1/(1\vee\upsilon (X))$, 
and they also calculated asymptotic bias and variance.

\paragraph*{Oscillation $\eta$-summing index of stochastic processes.}
Here we show that the conditions of Theorem \ref{main} hold for almost
all sample functions of several important classes of stochastic processes.
To this aim we use known results on Hausdorff-Besicovitch
dimension of graphs of sample functions.
Let $E\subset\RR^2$ be a bounded set, and let $\diam (A)$ denote the
diameter of a set $A\subset\RR^2$.
An $\epsilon$-covering of $E$ is a countable collection $\{E_k\colon\,
k\geq 1\}$ of sets such that $E\subset\cup_kE_k$ and $\sup_k\diam (E_k)
\leq\epsilon$.
For $s>0$, the Hausdorff $s$-measure of $E$ is defined by
$$
\CAH^s(E):=\lim_{\epsilon\downarrow 0}\inf\Big\{\sum_{k\geq 1}
\big (\diam (E_k)\big )^s\colon\, \{E_k\colon\,k\geq 1\}\,\,
\mbox{is an $\epsilon$-covering of $E$}\,\,\Big\}.
$$
Given $E\subset\RR^2$, the function $s\mapsto\CAH^s (E)$ is nonincreasing.
In fact, there is a critical value $s_c$ such that $\CAH^s (E)
=\infty$ for $s<s_c$ and $\CAH^s (E)=0$ for $s>s_c$.
This critical value $s_c$ is called the 
{\em Hausdorff-Besicovitch dimension} and is denoted by $\dim_{HB}(E)$. 
That is,
$$
\dim_{HB}(E):=\inf\big\{s>0\colon\,\CAH^s (E)=0\}
=\sup\{s>0\colon\,\CAH^s(E)=\infty\}.
$$
A relation between the lower metric entropy index $\Delta^{-}(E)$ of $E$
and the Hausdorff-Besicovitch dimension of $E$ is given by the following
result.

\begin{lem}\label{HB-BC}
For a bounded subset $E\subset\RR^2$, $\dim_{HB}(E)\leq\Delta^{-}(E)$.
\end{lem}

\begin{proof}
Let $E\subset\RR^2$ be bounded, and let $s>\Delta^{-}(E)$.
Since 
$$
\Delta^{-}(E)=\sup\big\{\alpha >0\colon\,\lim_{\epsilon\downarrow 0}
N(E;\epsilon)\epsilon^{\alpha}=+\infty\big\},
$$
we have that $\lim_{\epsilon\downarrow 0}N(E;\epsilon)\epsilon^{s}=0$.
Thus, for each $\epsilon >0$, there exists an  $\epsilon$-covering of $E$ by 
$N (E;\epsilon)$ balls of equal diameter $\epsilon$.
Hence
$$
\inf\Big\{\sum_{k\geq 1}\big (\diam (E_k)\big )^{s}\colon\, 
\{E_k\colon\,k\geq 1\}\,\,\mbox{is an $\epsilon$-covering of $E$}\,\,\Big\}
\leq N (E;\epsilon)\epsilon^{s},
$$
which gives $\CAH^{s}(E)\leq\liminf_{\epsilon\downarrow 0}
N (E;\epsilon)\epsilon^{s}=0$.
Thus $\dim_{HB}(E)\leq s$, and so $\dim_{HB}(E)\leq\Delta^{-}(E)$,
proving the lemma.
\qed\end{proof}

Let $\gr (X)$ be the graph of a regulated sample function of a stochastic 
process $X$.
By Proposition \ref{3main} and by the preceding lemma, we have
\beq\label{4main}
\dim_{HB}(\gr (X))\leq\Delta^{-}(\gr (X))\leq\Delta^{+}(\gr (X))
\leq 2-1/(1\vee\upsilon (X)).
\eeq
In fact, the left side is equal to the right side almost surely
for several classes of stochastic processes.
For example, let $X_{\alpha}$ be a symmetric $\alpha$-stable process for
some $0<\alpha\leq 2$.
Then by Theorem B of Blumenthal and Getoor (1962), for
almost all sample functions of $X_{\alpha}$, we have
$$
\dim_{HB}(\gr (X_{\alpha}))=2-\frac{1}{1\vee\alpha}
=2-\frac{1}{1\vee\upsilon (X_{\alpha})}.
$$
For another example,
let $X$ be a stochastic process of Example \ref{exmp2} having
the Orey index $\gamma_X\in (0,1)$.
Then by Theorem 1 of Orey (1970), for almost all sample functions of
$X$, we have
$$
\dim_{HB}(\gr (X))=2-\gamma_X=2-\frac{1}{\upsilon (X)}.
$$
For these stochastic processes, relation (\ref{4main}) yields that assumption
(\ref{2main}) of Theorem \ref{main} is satisfied for almost every
sample function, and so we have the following result. 

\begin{cor}\label{HBdim=var}
Let $X=\{X(t)\colon\,0\leq t\leq 1\}$ be a stochastic process, and let
$\eta=\{N_m\colon\,m\geq 1\}$ be a sequence of strictly increasing positive 
integers.
The relation
$$
\delta_{\eta}(X)\equiv
\lim_{m\to\infty}\frac{\log Q(X;\lambda (m))/N_m}{\log (1/N_m)}
=\frac{1}{1\vee\upsilon (X)}
$$
holds for almost all sample functions of $X$ provided 
either $(a)$ or $(b)$ holds, where
\begin{enumerate}
\item[$(a)$] $X$ is a symmetric $\alpha$-stable process for some
$\alpha\in (0,2]${\rm ;}
\item[$(b)$] $X$ is a mean zero Gaussian stochastic process with
stationary increments continuous in quadratic mean and such that
the Orey index $\gamma_X$ exists.
\end{enumerate}
\end{cor}

Similar results for more general processes other than $(a)$ and $(b)$
of Corollary \ref{HBdim=var} can be respectively found in
Pruitt and Taylor (1969, Section 8) and in K\^ono (1986). 

\section{Simulated symmetric $\alpha$-stable process}\label{simul:sas}

In this section we carry out a simulation study of small-sample
properties of the OS estimators from Definition \ref{ose}. 
To this aim, we simulate a symmetric $\alpha$-stable process for
several values of the exponent $\alpha$, which is equal to
its $p$-variation index.
Using repeated samples we calculate the bias, the standard deviation and 
the mean square error for the two estimators.

\paragraph*{Simulating $\sas$ process.}
Let $X_{\alpha}=\{X_{\alpha}(t)\colon\,t\geq 0\}$ be a symmetric
$\alpha$-stable stochastic process with the exponent $\alpha\in (0,2]$.
As stated in the introduction, the $p$-variation index of $X_{\alpha}$
is given by $\upsilon (X_{\alpha})=\alpha$ almost surely.
Since the OS estimators do not capture the values
of the $p$-variation index below $1$, we restrict our study to
estimating the exponent $\alpha \in [1,2]$.
To simulate a sample function of a $\sas$ process $X_{\alpha}$,
we generate a set $\{\xi_i\colon\,i=1,\dots,n\}$ of symmetric 
$\alpha$-stable pseudo-random variables and use the central limit theorem 
to get an approximation $\widetilde{X}_{\alpha}$ of $X_{\alpha}$, where
$$
\widetilde{X}_{\alpha}(t)=\frac{1}{n^{1/\alpha}}\sum_{i=1}^{[nt]} \xi_i,
\qquad 0\leq t\leq 1,
$$
and $[r]$ denotes the integer part of $r$.
By the central limit theorem, the distribution of 
$\widetilde{X}_{\alpha}$ on the Skorohod space $D[0,1]$ converges
weakly to the distribution of $X_{\alpha}$ as $n\to\infty$.
We take $n=2^{14}$.
To generate a symmetric $\alpha$-stable random variable $\xi$, we use the
results of Chambers et al. (1976) (see also Section 4.6
in Zolotarev 1986).
That is,  in the sense of equality in distribution, we have
$$
\xi=\frac{\sin\alpha U}{(\cos U)^{1/\alpha}}\cdot
\Big (\frac{\cos (U-\alpha U)}{E}\Big )^{(1-\alpha )/\alpha}
\quad\mbox{if $\alpha\not =1$,}\qquad\mbox{and}\qquad \xi =\tan U\quad
\mbox{if $\alpha =1$,}
$$
where the random variable $U$ has the uniform distribution on 
$[-\pi/2,\pi/2]$, and the random variable $E$, which is independent of $U$, 
has the standard exponential distribution.
All calculations are done using the computing system {\em Mathematica}.

\begin{table}[tb]
\caption{Properties of $100$ samples of estimates $\naive (X_{\alpha})$
based on $\{2^m\colon\,1\leq m\leq 14\}$}
\label{simul1}
\medskip
\centering
\small{\begin{tabular}{|c||c|c|c|c|c|c|c|c|c|c|c|}\hline
$\alpha$& 1.0 & 1.1 & 1.2 & 1.3 & 1.4 & 1.5 & 1.6 & 1.7 & 1.8 & 1.9 & 2.0\\ 
\hline\hline 
\rule[0mm]{0mm}{4mm} $\overline{\alpha}$ & 1.278  & 1.346 & 1.406 & 1.473
        & 1.556  & 1.632 & 1.715  & 1.798 & 1.883 & 1.967 & 2.051 \\ \hline
\rule[0mm]{0mm}{4mm} $\overline{\alpha}-\alpha$ & .2783 & .2464 & .2064 & .1727
        & .1563 & .1322  & .1149  & .0984 & .0832 & .0667 & .0510\\ \hline
\rule[0mm]{0mm}{4mm} $SD$ & .0867 & .0864 & .0542 & .0285 & .0264
        & .0147 & .0089 & .0063 & .0059 & .0040 & .0025 \\ \hline
\rule[0mm]{0mm}{4mm} $MSE$ & .0849 & .0681 & .0455 & .0306 & .0251        
        & .0177 & .0133 & .0097 & .0070 & .0045 & .0026 \\ \hline
\end{tabular}}
\end{table}

\begin{table}[tb]
\caption{Properties of $100$ samples of estimates $\ols (X_{\alpha})$
based on $\{2^m\colon\,1\leq m\leq 14\}$}
\label{simul2}
\medskip
\centering
\small{\begin{tabular}{|c||c|c|c|c|c|c|c|c|c|c|c|}\hline
$\alpha$& 1.0 & 1.1 & 1.2 & 1.3 & 1.4 & 1.5 & 1.6 & 1.7 & 1.8 & 1.9 & 2.0\\ 
\hline\hline 
\rule[0mm]{0mm}{4mm} $\overline{\alpha}$ & 1.158  & 1.197 & 1.269 & 1.337
      & 1.383  & 1.464 & 1.529  & 1.598 & 1.654 & 1.711 & 1.767 \\ \hline
\rule[0mm]{0mm}{4mm} $\overline{\alpha}-\alpha$ & .1582 & .0974 & .0693 & .0375
      & -.0172 & -.0362  & -.0712  & -.1022 & -.1459 & -.1887 & -.2329\\ \hline
\rule[0mm]{0mm}{4mm} $SD$ & .0726 & .0785 & .0902 & .0804 & .0964
      & .0794 & .0791 & .0786 & .0691 & .0555 & .0435 \\ \hline
\rule[0mm]{0mm}{4mm} $MSE$ & .0302 & .0156 & .0129 & .0078 & .0095
      & .0076 & .0113 & .0166 & .0260 & .0387 & .0560 \\ \hline
\end{tabular}}
\end{table}

\paragraph*{OS estimators.}
Let $\eta=\{2^m\colon\,m\geq 1\}$.
We simulate a sample function $\widetilde{X}_{\alpha}$ at $2^{14}+1$ equally 
spaced points  $\{u_1,\dots,u_N\}=\{i2^{-14}\colon\,i=0,\dots,2^{14}\}$.
For each $1\leq m\leq 14$, let
$$
Q(m):=\sum_{i=1}^{2^m}\Big [
\max_{u_k\in\Delta_{i,m}}\big\{\widetilde{X}_{\alpha}(u_k)\big\}
-\min_{u_k\in\Delta_{i,m}}\big\{\widetilde{X}_{\alpha}(u_k)\big\}\Big ],
$$
where $\Delta_{i,m}:=[(i-1)2^{-m},i2^{-m}]$. 
Then define $r(m)$, $m=1,\dots,14$, by relation (\ref{r_m}) with $N_m=2^m$.
We use the OS estimators based on $\eta_{14}=\{2^m\colon\,1\leq m\leq 14\}$.  
Thus by Definition \ref{ose}, the naive oscillation $\eta$-summing estimator 
$\naive (X_{\alpha})=r(14)$, and
the OLS oscillation $\eta$-summing estimator
$$
\ols (X_{\alpha})=\frac{\sum_{m=1}^{14}(x_m-\bar x)m}{
\sum_{m=1}^{14}(x_m-\bar x)^2},
$$
where $x_m=\log_2(2^m/Q(m))$ and $\bar x=14^{-1}\sum_{m=1}^{14}x_m$.

\paragraph*{Monte-Carlo study.}
For $11$ different values of $\alpha\in\{1.0,1.1,\dots,2.0\}$, 
we simulate a vector 
$\{\widetilde{X}_{\alpha}(k/2^M): k=0,\dots,2^M\}$ of values
of a sample function of $X_{\alpha}$ with $M = 14$.
By the preceding paragraph, we have two estimators $\widehat{\alpha}$ of 
$\alpha$: the naive oscillation $\eta$-summing estimator $\naive (X_{\alpha})$,
and  the OLS oscillation $\eta$-summing estimator $\ols (X_{\alpha})$,
both based on $\{2^m\colon\,1\leq m\leq 14\}$.   
We repeat this procedure $K=100$ times to obtain the estimates 
$\widehat{\alpha}_1,\dots \widehat{\alpha}_K$ of $\alpha$ 
for each of the two cases. 
Then we calculate:
\begin{itemize}
\item The estimated expected value $\overline{\alpha}
      :=\big (\sum_{i=1}^K\widehat{\alpha}_i\big )/K$;
\item The bias $\overline{\alpha}-\alpha$;
\item The estimated standard deviation  $SD :=\sqrt{\sum_{i=1}^K
      (\widehat{\alpha}_i-\overline{\alpha})^2/(K-1)}$;
\item The estimated mean square error $MSE :=
\sum_{i=1}^K(\widehat{\alpha}_i-\alpha)^2/K$.
\end{itemize}
The estimation results are presented in Tables \ref{simul1} and \ref{simul2}.  
Next is a qualitative description of the performance of the two OS
estimators.

\begin{description}
\item[Bias] The estimation results show a different behavior of the bias
for the two estimators.
The naive estimator $\naive(X_{\alpha})$ display monotonically decreasing 
positive bias when $\alpha$ values increase from 1 to 2.
While the bias of the OLS estimator $\ols (X_{\alpha})$ monotonically 
decrease from a positive bias for $\alpha < 1.4$ to a negative bias 
for $\alpha\geq 1.4$, and the minimal absolute bias is achieved
when $\alpha =1.4$.
\end{description}
\noindent
\begin{description}
\item[SD] The estimated standard deviation is also different for
the two estimators.
SD values for the naive estimator $\naive (X_{\alpha})$ monotonically
decrease when $\alpha$ values increase from 1 to 2, and are quite
small when $\alpha$ is close to $2$.
While SD values for the OLS estimator $\ols (X_{\alpha})$ remain similar
and somewhat larger than for the naive estimator; only a little 
improvement one can notice when $\alpha$ is close to 2.
\end{description}
\noindent
\begin{description}
\item[MSE] The estimated mean square error remain different for
the two estimators.
MSE values for the naive estimator $\naive (X_{\alpha})$ decrease steady
when $\alpha$ values increase from 1 to 2, 
while MSE values for the OLS estimator $\naive (X_{\alpha})$ is smallest 
when $\alpha =1.5$,
and are increasing for all other values of $\alpha$.
\end{description}
\noindent
In conclusion the results show a distinction between the two OS estimators:
the OLS estimator $\naive (X_{\alpha})$ display better
performance for values $1.4\leq\alpha\leq 1.6$,
while the naive estimator $\naive (X_{\alpha})$
behaves best for $\alpha$ values close to  $2$.
This Monte-Carlo study was extended to sample functions based on
a larger number of points, 
i.e.\ $\{\widetilde{X}_{\alpha}(k/2^M): k=0,\dots,2^M\}$ with 
$M\in\{15,16,17\}$.
The results from the increased sample size show the same qualitative 
behavior as before, but with increased accuracy 
(see Norvai\v sa and Salopek, 2000b).

To the best of our knowledge, the two OS estimators provide the first
attempt to estimate the exponent $\alpha$ of a $\sas$ process from
a sample function.
Recently Crovella and Taqqu (1999) introduced a new method to estimate 
the exponent $\alpha$ of a $\sas$ random variable.

\section{Simulated fractional Brownian motion}\label{simul:fbm}

Here we perform a Monte Carlo study like in the preceding section,  
with a fractional Brownian motion.
A fractional Brownian motion $\{B_H(t)\colon\,t\geq 0\}$ with the Hurst 
exponent $H\in (0,1)$ is a Gaussian stochastic process with stationary 
increments having the covariance function
$$
E\big [B_H(t)B_H(s)\big ]=\frac{1}{2}\Big [t^{2H}+s^{2H}-|t-s|^{2H}\Big ],
\qquad\mbox{for $t,s\geq 0$}
$$
and $B_H(0)=0$ almost surely.
Since its incremental variance $\sigma_{B_H}(u)=u^H$, the Orey index
$\gamma_{B_H}$ exists and is equal to the Hurst exponent $H$
(cf. Example \ref{exmp2}). 
Thus by relation (\ref{var-of-Orey}), $B_H$ has the $p$-variation index 
$\upsilon (B_H)=1/H$ almost surely.
To simulate a fractional Brownian motion we use the program of
Maeder (1995) written in {\em Mathematica}.

\begin{table}[tb]
\caption{Properties of $100$ samples of estimates $\naive (B_{H})$
based on $\{2^m\colon\,1\leq m\leq 14\}$}
\label{simul5}
\medskip
\centering
\small{\begin{tabular}{|c||c|c|c|c|c|c|c|c|c|c|}\hline
$h=1/H$& 1.2 & 1.4 & 1.6 & 1.8 & 2.0 & 2.2 & 2.4 & 2.6 & 2.8 & 3.0 \\ 
\hline\hline 
\rule[0mm]{0mm}{4mm} $\overline{h}$ & 1.132  & 1.304 & 1.474 & 1.640
     & 1.802  & 1.961 & 2.117  & 2.270 & 2.420 & 2.566 \\ \hline
\rule[0mm]{0mm}{4mm} $\overline{h}-h$ & -.0683 & -.0956 & -.1262 & -.1603
     & -.1979 & -.2390  & -.2826  & -.3296 & -.3804 & -.4339 \\ \hline
\rule[0mm]{0mm}{4mm} $SD$ & .0039 & .0018 & .0013 & .0020 & .0021
        & .0027 & .0026 & .0035 & .0042 & .0042 \\ \hline
\rule[0mm]{0mm}{4mm} $MSE$ & .0047 & .0092 & .0159 & .0257 & .0392        
        & .0571 & .0799 & .1086 & .1447 & .1883 \\ \hline
\end{tabular}}
\end{table}

\begin{table}[tb]
\caption{Properties of $100$ samples of estimates $\ols (B_{H})$
based on $\{2^m\colon\,1\leq m\leq 14\}$}
\label{simul6}
\medskip
\centering
\small{\begin{tabular}{|c||c|c|c|c|c|c|c|c|c|c|}\hline
$h=1/H$& 1.2 & 1.4 & 1.6 & 1.8 & 2.0 & 2.2 & 2.4 & 2.6 & 2.8 & 3.0 \\ 
\hline\hline 
\rule[0mm]{0mm}{4mm} $\overline{h}$ & 1.186  & 1.340 & 1.481 & 1.614
    & 1.712  & 1.821 & 1.891  & 1.971 & 2.043 & 2.095 \\ \hline
\rule[0mm]{0mm}{4mm} $\overline{h}-h$ & -.0143 & -.0598 & -.1185 & -.1864
    & -.2879 & -.3795  & -.5089  & -.6294 & -.7572 & -.9053 \\ \hline
\rule[0mm]{0mm}{4mm} $SD$ & .0490 & .0466 & .0499 & .0469 & .0465
        & .0474 & .0399 & .0447 & .0494 & .0471 \\ \hline
\rule[0mm]{0mm}{4mm} $MSE$ & .0026 & .0057 & .0165 & .0369 & .0850        
        & .1462 & .2606 & .3982 & .5757 & .8217 \\ \hline
\end{tabular}}
\end{table}

\paragraph{OS and G estimation.}
In this section we apply four estimators of the $p$-variation index.
As before, the two OS estimators from Definition \ref{ose} will be used 
to estimate $h:=1/H$.
Moreover, we invoke the two estimators of the Orey index
introduced in Norvai\v sa and Salopek (2000), and which will be called
{\em G estimation}, which is short for the Gladyshev estimation.
More specifically, let $\eta=\{2^m\colon\,m\geq 1\}$.
We simulate a sample function $\widetilde{B}_H$ at $2^{14}+1$ equally
spaced points $\{u_1,\dots,u_N\}=\{i2^{-14}\colon\,i=0,\dots,2^{14}\}$.
For each $1\leq m\leq 14$, let
$$
Q(m):=\sum_{i=1}^{2^m}\Big [
\max_{u_k\in\Delta_{i,m}}\big\{\widetilde{B}_H(u_k)\big\}
-\min_{u_k\in\Delta_{i,m}}\big\{\widetilde{B}_H(u_k)\big\}\Big ],
$$
where $\Delta_{i,m}:=[(i-1)2^{-m},i2^{-m}]$. 
Then define $r(m)$, $m=1,\dots,14$, by (\ref{r_m}) with $N_m=2^m$.
Therefore one can use estimators based on $\eta_{14}=\{2^m\colon\,
1\leq m\leq 14\}$.  
Thus by Definition \ref{ose}, the naive oscillation $\eta$-summing estimator 
$\naive (B_H)=r(14)$, and the OLS oscillation $\eta$-summing estimator 
$$
\ols (B_H)=\frac{\sum_{m=1}^{14}(x_m-\bar x)m}{
\sum_{m=1}^{14}(x_m-\bar x)^2},
$$
where $x_m=\log_2(2^m/Q(m))$ and $\bar x=14^{-1}\sum_{m=1}^{14}x_m$.

To recall the G estimation, again let $\eta=\{2^m\colon\,m\geq 1\}$, and
let $\widetilde{B}_H$ be a sample function given by its values at
$2^{14}+1$ equally spaced points $\{i2^{-14}\colon\,i=0,\dots,2^{14}\}$.
For each $1\leq m\leq 14$, let
$$
s_2(m):=\sum_{i=1}^{2^m}\Big [\widetilde{B}_H(i/N_m)
- \widetilde{B}_H((i-1)/N_m)\Big ]^2\quad\mbox{and}\quad
r(m):=\frac{\log 2^{-m}}{\log \sqrt{s_2(m)2^{-m}}}.
$$
The {\em naive Gladyshev estimator} of the $p$-variation index
is defined by $\gnaive (B_H)=r(14)$, and  the OLS Gladyshev estimator 
is defined by
$$
\gols (B_H)=\frac{\sum_{m=1}^{14}(x_m-\bar x)m}{
\sum_{m=1}^{14}(x_m-\bar x)^2},
$$
where $x_m=\log_2\sqrt{2^m/s_2(m)}$ and $\bar x=14^{-1}\sum_{m=1}^{14}x_m$.

\begin{table}[tb]
\caption{Properties of $100$ samples of estimates $\gnaive (B_{H})$
based on $\{2^m\colon\,1\leq m\leq 14\}$}
\label{simul7}
\medskip
\centering
\small{\begin{tabular}{|c||c|c|c|c|c|c|c|c|c|c|}\hline
$h=1/H$& 1.2 & 1.4 & 1.6 & 1.8 & 2.0 & 2.2 & 2.4 & 2.6 & 2.8 & 3.0 \\ 
\hline\hline 
\rule[0mm]{0mm}{4mm} $\overline{h}$ & 1.164  & 1.348 & 1.530 & 1.709
    & 1.887  & 2.062 & 2.236  & 2.407 & 2.576 & 2.742 \\ \hline
\rule[0mm]{0mm}{4mm} $\overline{h}-h$ & -.0364 & -.0524 & -.0704 & -.0907
    & -.1133 & -.1383  & -.1645  & -.1928 & -.2245 & -.2580 \\ \hline
\rule[0mm]{0mm}{4mm} $SD$ & .0039 & .0020 & .0013 & .0023 & .0022
        & .0031 & .0030 & .0039 & .0044 & .0048 \\ \hline
\rule[0mm]{0mm}{4mm} $MSE$ & .0013 & .0027 & .0050 & .0082 & .0128        
        & .0191 & .0271 & .0372 & .0504 & .0666 \\ \hline
\end{tabular}}
\end{table}

\begin{table}[tb]
\caption{Properties of $100$ samples of estimates $\gols (B_{H})$
based on $\{2^m\colon\,1\leq m\leq 14\}$}
\label{simul8}
\medskip
\centering
\small{\begin{tabular}{|c||c|c|c|c|c|c|c|c|c|c|}\hline
$h=1/H$& 1.2 & 1.4 & 1.6 & 1.8 & 2.0 & 2.2 & 2.4 & 2.6 & 2.8 & 3.0 \\ 
\hline\hline 
\rule[0mm]{0mm}{4mm} $\overline{h}$ & 1.236  & 1.429 & 1.628 & 1.818
     & 1.987  & 2.184 & 2.325  & 2.483 & 2.645 & 2.790 \\ \hline
\rule[0mm]{0mm}{4mm} $\overline{h}-h$ & .0363 & .0292 & .0281 & .0185
     & -.0131 & -.0161  & -.0751  & -.1171 & -.1548 & -.2100 \\ \hline
\rule[0mm]{0mm}{4mm} $SD$ & .0705 & .0870 & .1209 & .1241 & .1418
        & .1689 & .1425 & .2459 & .2650 & .2343 \\ \hline
\rule[0mm]{0mm}{4mm} $MSE$ & .0062 & .0084 & .0153 & .0156 & .0201        
        & .0285 & .0258 & .0736 & .0935 & .0985 \\ \hline
\end{tabular}}
\end{table}

\paragraph{Monte-Carlo study.}
For $10$ different values of $H\in\{0.83\approx 1.2^{-1},0.71\approx 
1.4^{-1},\dots,0.33\approx 3.0^{-1}\}$,
we simulate the vector $\{ B_H(i2^{-m}): i=0,\dots,2^m\}$ with $m = 14$,
and calculate the four estimators.
This gives us four different estimates $\widehat{h}$ of $h=1/H$.  
We repeat this procedure $K=100$ times to obtain the estimates 
$\widehat{h}_1,\dots \widehat{h}_K$ of $h$ for each of the four cases. 
Then we calculate:
\begin{itemize}
\item The estimated expected value $\overline{h}
      :=\big (\sum_{i=1}^K\widehat{h}_i\big )/K$;
\item The bias $\overline{h}-h$;
\item The estimated standard deviation  $SD :=\sqrt{\sum_{i=1}^K
      (\widehat{h}_i-\overline{h})^2/(K-1)}$;
\item The estimated mean square error $MSE :=
\sum_{i=1}^K(\widehat{h}_i-h)^2/K$.
\end{itemize}
The estimation results are presented in Tables \ref{simul5} - \ref{simul8}.  
First one can compare the OS estimation results of $\alpha\in [1,2]$ 
(Tables \ref{simul1} and \ref{simul2}) for a $\sas$ process, 
and the OS estimation results of
$h=1/H\in (1,2]$ (columns $h=1.2,1.4,1.6,1.8,2.0$ of Tables 
\ref{simul5} and \ref{simul6}) for a fractional Brownian motion with
the Hurst exponent $H$.
The accuracy of the naive OS estimator for the two processes is similar.
There is only some differences in the character of monotonicity along
different values of parameters.  
The same holds for the OLS OS estimator for the two processes.

Now if we look at columns $h\in\{2.0,2.2,\dots,3.0\}$ of Tables \ref{simul5} 
and \ref{simul6}, it is clear that the OS estimates of $h$ 
are very poor in this case.
This is so since estimation errors appear in the denominator of
the relation (\ref{r_m}).
For example, if $H$ and $\widehat{H}$ both are small, then the left side of
the relation
\beq\label{difference}
\frac{1}{H}-\frac{1}{\widehat{H}}=(\widehat{H}-H)\frac{1}{H\widehat{H}}
\eeq
can be relatively large as compared to $\widehat{H}-H$.
The same remark applies to the G estimates in Tables \ref{simul7} and 
\ref{simul8}. 
To show that this is so we applied the OS estimator to evaluate the Hurst
exponent directly using the relations (\ref{var-of-Orey}) and 
(\ref{index=var}).
This means that in our earlier estimation formulas we need just
to interchange the numerator and the denominator in (\ref{r_m}).
The estimation results are presented in Tables \ref{simul3} and 
\ref{simul4}, where $\widehat{H}$, the bias, SD and MSE are defined as
before with $h$ replaced by $H$.
In Tables \ref{simul5} - \ref{simul4}, the estimation results for
$H=1/2$, and so for $h=2$, are all based on the same set of $100$ simulated
sample functions, which can be used to verify the effect of the above 
relation \ref{difference}.
Also, the results of the OS estimation of $H$ in Tables 
\ref{simul3} and \ref{simul4} can be compared with the results of the
G estimation of $H$ in Tables 1 and 3 of Norvai\v sa and Salopek (2000).
The naive estimators corresponding to the OS and G estimations
show very similar properties.
As far as the OLS estimators concern, for small vaules of $H$,
the bias and MSE of the OS estimation are larger than the bias 
and MSE of G estimation.
However, for the same values of $H$, the standard deviation of
the OS estimation is smaller than the standard deviation of the
G estimation. 
In sum the two estimation methods OS and G show similar results
when applied to estimate the Hurst exponent of a fractional
Brownian motion.

\begin{table}[tb]
\caption{Properties of $100$ samples of estimates $\hnaive (B_{H})$
based on $\{2^m\colon\,1\leq m\leq 14\}$}
\label{simul3}
\medskip
\centering
\small{\begin{tabular}{|c||c|c|c|c|c|c|c|c|c|}\hline
$H$& 0.1 & 0.2 & 0.3 & 0.4 & 0.5 & 0.6 & 0.7 & 0.8 & 0.9  \\ 
\hline\hline 
\rule[0mm]{0mm}{4mm} $\overline{H}$ & .1578  & .2572 & .3566 & .4559
     & .5549  & .6539 & .7526  & .8511 & .9496  \\ \hline
\rule[0mm]{0mm}{4mm} $\overline{H}-H$ & .0578 & .0572 & .0566 & .0559
     &  .0549 & .0539  & .0526  & .0511 & .0496  \\ \hline
\rule[0mm]{0mm}{4mm} $SD$ & .0006 & .0007 & .0006 & .0007 & .0006
        & .0007 & .0010 & .0020 & .0101 \\ \hline
\rule[0mm]{0mm}{4mm} $MSE$ & .0033 & .0033 & .0032 & .0231 & .0030        
        & .0029 & .0027 & .0026 & .0026  \\ \hline
\end{tabular}}
\end{table}

\begin{table}[tb]
\caption{Properties of $100$ samples of estimates $\hols (B_{H})$
based on $\{2^m\colon\,1\leq m\leq 14\}$}
\label{simul4}
\medskip
\centering
\small{\begin{tabular}{|c||c|c|c|c|c|c|c|c|c|}\hline
$H$& 0.1 & 0.2 & 0.3 & 0.4 & 0.5 & 0.6 & 0.7 & 0.8 & 0.9  \\ 
\hline\hline 
\rule[0mm]{0mm}{4mm} $\overline{H}$ & .3471  & .3932 & .4481 & .5109
    & .5802  & .6508 & .7320  & .8146 & .8962 \\ \hline
\rule[0mm]{0mm}{4mm} $\overline{H}-H$ & .2471 & .1932 & .1481 & .1109
    & .0802 & .0508  & .0320  & .0146 & -.0038 \\ \hline
\rule[0mm]{0mm}{4mm} $SD$ & .0069 & .0081 & .0113 & .0141 & .0175
        & .0235 & .0276 & .0364 & .0449 \\ \hline
\rule[0mm]{0mm}{4mm} $MSE$ & .0611 & .0374 & .0221 & .0125 & .0067        
        & .0031 & .0018 & .0015 & .0020 \\ \hline
\end{tabular}}
\end{table}

\section{Financial data analysis}

In this section, we analyze the financial data set provided by 
Olsen \& Associates.
It is the high-frequency data set HFDF96, which consists of $25$ different
foreign exchange spot rates, $4$ spot metal rates, and $2$ series of
stock indices.
This data set was recorded from $1$ Jan 1996 GMT to 31 Dec 1996 GMT.
Each set has 17568 entries recorded at half hour intervals.
The same financial data set was studied in Norvai\v sa and Salopek (2000) 
using the two estimators of the Orey index based on the result
of Gladyshev (1961).

\paragraph*{Returns.}
First notice that returns in continuous time and 
discrete time financial models are treated slightly differently.
A return in a discrete time model is a function 
$\widehat{R}$ defined on a lattice $t\in\{0,1,\dots,T\}$ with values 
being a suitable transform of a pair $\{P(t-1),P(t)\}$, where $P$ is a 
stock price process.
A return in a continuous time model is a function $R$ defined
on $[0,T]$ so that
$R(t)-R(t-1)=\widehat{R}(t)$ and $R(0)=0$ for all $t\in \{1,\dots,T\}$.
This gives a $1-1$ correspondence between continuous time returns
used in this paper and the usual discrete time returns
(see Section 2.1 in Norvai\v sa, 2000a, for further details).

Given a historical data set $\{d_0,\dots,d_K\}$ of values of a financial 
asset, let $P$ be a function defined on $[0,1]$ with values 
$d_k$ at $u_k:=k/K$ for $k=0,\dots,K$.
Usually in econometric literature, returns are log transforms of
the price.
Thus in continuous time models, this corresponds to assuming that
$P(t)=P(0)\exp\{X(t)\}$, $0\leq t\leq 1$, for some return process $X$
which is to be analyzed.
Alternatively, the price process $P$ can be a solution of an integral
equation, which is not a simple exponential.
If a stochastic process $X$ has the quadratic variation along the sequence
of partitions $\lambda=\{\lambda (m)\colon\,m\geq 1\}$ defined by 
(\ref{1ose}), and the price process $P$ is a solution of a linear integral 
equation with respect to $X$ 
(such as the Black-Scholes-Merton model (\ref{BSM-model}), 
with $B$ replaced by $X$), then $P$ also has the quadratic variation 
along the sequence $\lambda$ and
the process $X$ can be recovered by the relation 
$$
X(t)=R_{net}(P)(t):=\lim_{m\to\infty}\sum_{i=1}^{N_m}\big [P(t\wedge
i/N_m)-P(t\wedge(i-1)/N_m)\big ]/P(t\wedge (i-1)/N_m),
$$
for $0\leq t\leq 1$.
Here the quadratic variation is understood in the sense of F\"olmer (1981),
which is further developed in Norvai\v sa (2000b).
The above process $X$ will be called the {\em net return} of $P$,
which is analogous to discrete time simple net returns.
The {\em log return} $X$ of the price process $P$ is defined by 
$X(t)=R_{log}(P)(t):=\log [P(t)/P(0)]$, $0\leq t\leq 1$. 
The difference between the two returns is
\beq\label{dif_returns}
R_{net}(P)(t)-R_{log}(P)(t)=\frac{1}{2}\int_0^t\frac{d[P]^c}{P^2}
-\sum_{(0,t]}\Big [\log\frac{P}{P_{-}}-\frac{\Delta^{-}P}
{P_{-}}\Big ]-\sum_{[0,t)}\Big [\log\frac{P_{+}}{P}-\frac{\Delta^{+}P}{P}
\Big ]
\eeq
for $0\leq t\leq 1$, where $[P]^c$ is a continuous part of the quadratic 
variation of $P$.
Since $P$ has the quadratic variation along the sequence $\lambda$, 
the difference  $R_{net}(P)-R_{log}(P)$ given by (\ref{dif_returns}) has 
bounded variation. 
Thus the $p$-variation indices of the two returns are equal provided
both are not less than $1$.
A discussion of the difference $R_{net}(P)-R_{log}(P)$ when $P$ is the
geometric Brownian motion or $P$ is a model for USD/JPY exchange rates 
can be found on pages 362 and 366 of Norvai\v sa and Salopek (2000).
\begin{figure}
\onelinecaptionsfalse
\unitlength1cm
\centering
\begin{minipage}[b]{8.0cm}
 \centering
 \includegraphics[height=5.0cm, width=6.cm]{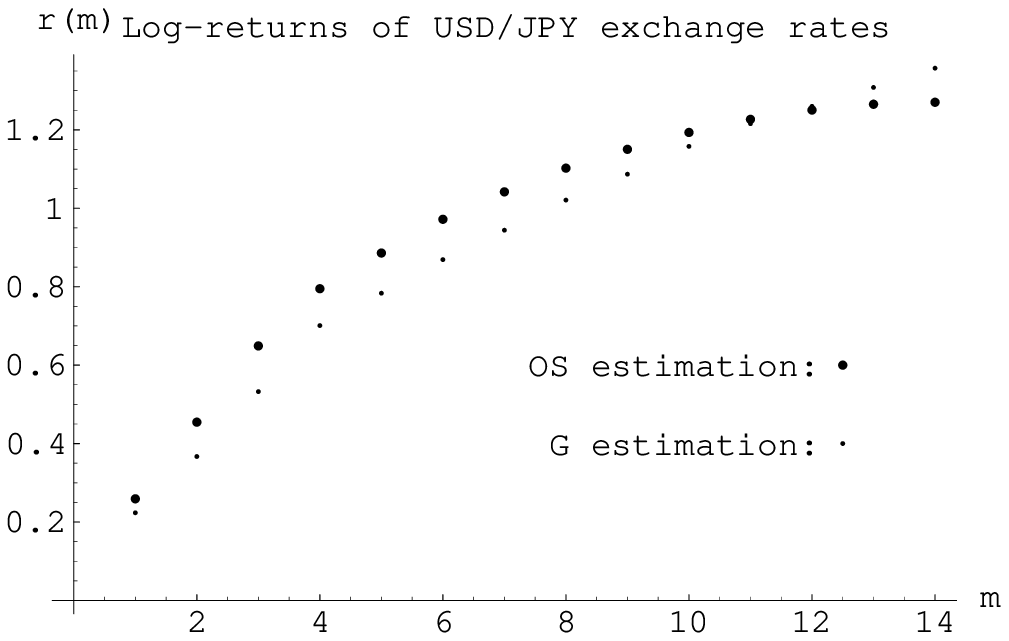}
\end{minipage}%
\hfill
\begin{minipage}[b]{8.0cm}
 \centering
 \includegraphics[height=5.0cm, width=6.cm]{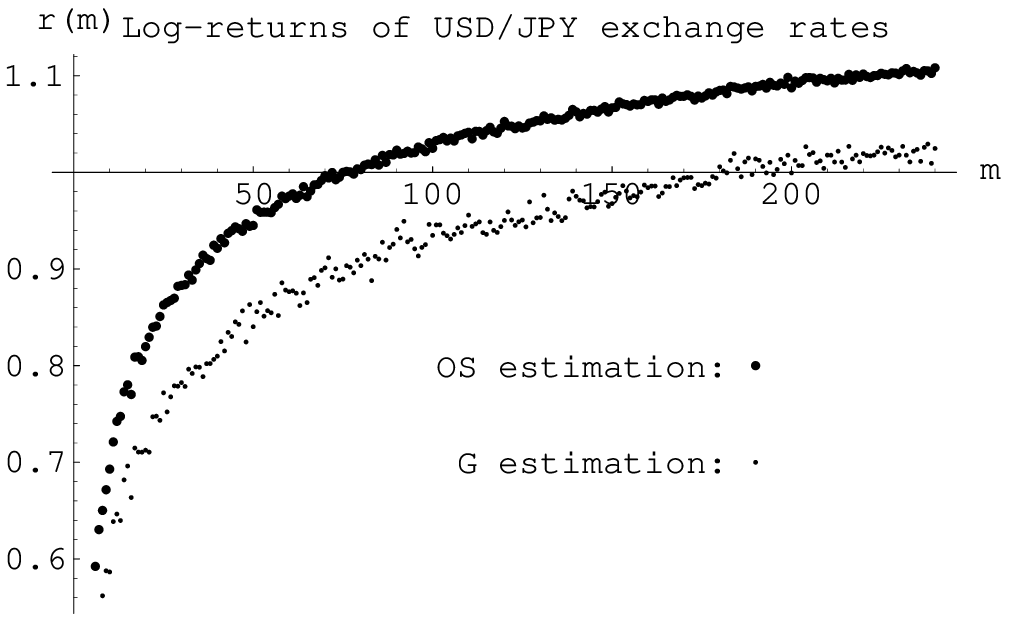}
\end{minipage}\\ 
\begin{minipage}[t]{8.0cm}
\captionstyle{flushleft}
\setcaptionwidth{8.cm}
\caption{Naive estimators: $\naive=1.27$ and $\gnaive=1.357$ 
based on $\{2^m\colon\,m=1,\dots,14\}$.}
\label{usdjpy1}
\end{minipage}% 
\hfill
\begin{minipage}[t]{8.0cm}
\setcaptionwidth{8.cm}
\caption{Naive estimators: $\naive=1.108$ and $\gnaive=1.025$
based on $\{m\colon\,m=1,\dots,240\}$.}
\label{usdjpy2}
\end{minipage}%
\end{figure}
\begin{figure}
\onelinecaptionsfalse
\unitlength1cm
\centering
\begin{minipage}[b]{8.0cm}
 \centering
 \includegraphics[height=5.0cm, width=6.cm]{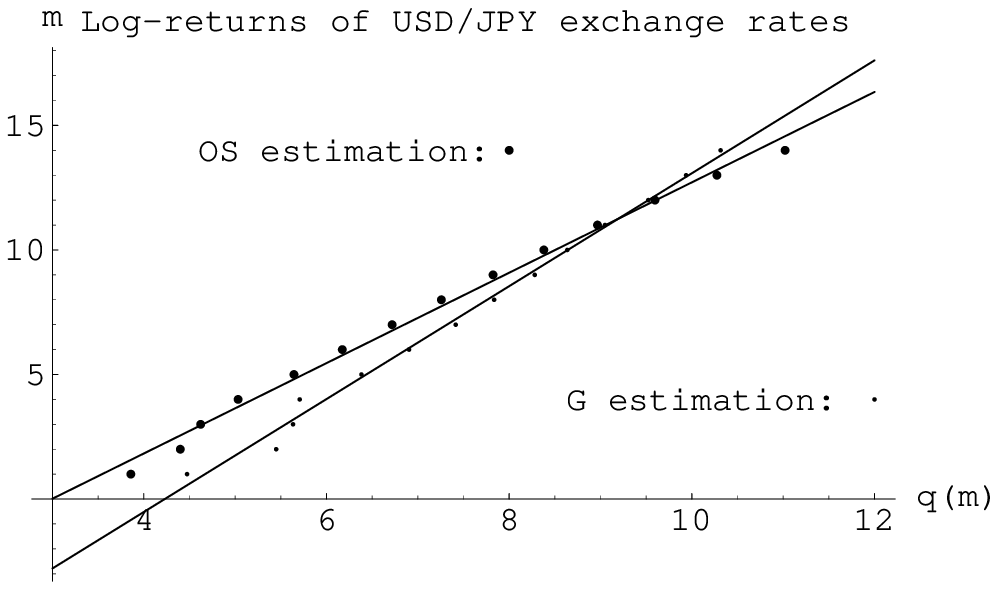}
\end{minipage}%
\hfill
\begin{minipage}[b]{8.0cm}
 \centering
 \includegraphics[height=5.0cm, width=6.cm]{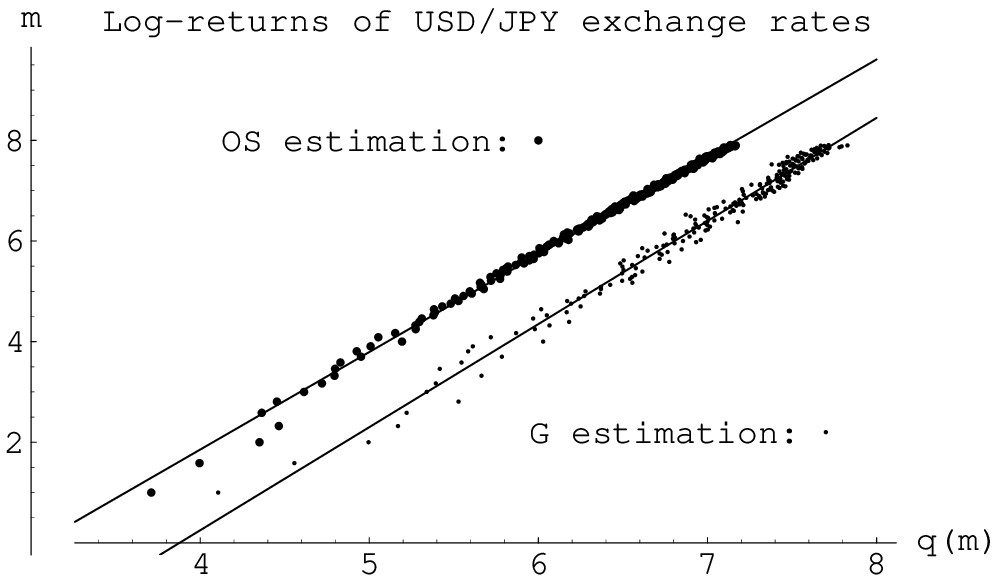}
\end{minipage}\\ 
\begin{minipage}[t]{8.0cm}
\captionstyle{flushleft}
\setcaptionwidth{8.cm}
\caption{OLS estimators: $\ols=1.814$ and $\gols=2.266$
based on $\{2^m\colon\,m=1,\dots,14\}$.}
\label{usdjpy3}
\end{minipage}% 
\hfill
\begin{minipage}[t]{8.0cm}
\setcaptionwidth{8.cm}
\caption{OLS estimators: $\ols =1.939$ and $\gols=2.049$
based on $\{m\colon\,m=1,\dots,240\}$.}
\label{usdjpy4}
\end{minipage}%
\end{figure}
\paragraph*{OS and G estimation.}
To estimate the $p$-variation index of a return, the oscillation
$\eta$-summing estimators from Definition \ref{ose} will be used
when the sequence $\eta=\{N_m\colon\,m\geq 1\}$ is given by
$N_m=2^m$ and $N_m=m$ for integers $m\geq 1$.
A comparison will be made with the results of Norvai\v sa and Salopek (2000)
by using their estimators based on the result of Gladyshev (1961).
As before, the oscillation $\eta$-summing estimation is called the
{\em OS estimation}, and the estimation as in Norvai\v sa and 
Salopek (2000) is called the {\em G estimation}.

Given a price process $P$ at points $u_k=k/K$, $k=0,\dots,K$,
find the maximal integer $M$ such that
\beq\label{1est}
\cup_{m=1}^M\lambda (m)\subset\{u_0,\dots,u_K\},\qquad
\mbox{where $\lambda (m)=\{i/N_m\colon\,i=0,\dots,N_m\}$.}
\eeq
Notice that in the case $N_m=m$, the sequence $\{\lambda (m)\colon\,
m\geq 1\}$ is {\em not} nested.
Since for HFDF96 data set $K=17568$, $M=14$ when $N_m=2^m$ and 
$M=240$ when $N_m=m$.

Let $\eta=\{N_m\colon\,m\geq 1\}$ be a sequence of strictly increasing 
positive integers, and let $M$ be such that relation (\ref{1est}) holds.
In the case $X=R_{net}(P)$ or $X=R_{log}(P)$,
for each $1\leq m\leq M$, let
$$
Q(m):=Q(X;\lambda (m))=\sum_{i=1}^{N_m}\Big [
\max_{u_k\in\Delta_{i,m}}\big\{X(u_k)\big\}
-\min_{u_k\in\Delta_{i,m}}\big\{X(u_k)\big\}\Big ],
$$
where $\Delta_{i,m}:=[(i-1)/N_m,i/N_m]$. Also for each $1\leq m\leq M$, let
\beq\label{2est}
q(m):=\log_2(N_m/Q(m))\quad\mbox{and}\quad
r(m):=\frac{\log (1/N_m)}{\log Q(m)/N_m}
=\frac{\log_2 N_m}{q(m)}.
\eeq
By Definition \ref{ose}, the naive oscillation $\eta$-summing estimator 
$\naive:=\naive (X)=r(M)$.
The OLS oscillation $\eta$-summing estimator $\ols$ is defined by
\beq\label{3est}
\ols:=\ols (X)=\frac{\sum_{m=1}^M(x_m-\bar x)
\log_2 N_m}{\sum_{m=1}^M(x_m-\bar x)^2},
\eeq
where $x_m=q(m)=\log_2(N_m/Q(m))$ and $\bar x=M^{-1}\sum_{m=1}^Mx_m$.

Turning to the G estimation, again let $\eta=\{N_m\colon\,m\geq 1\}$ 
be a sequence of strictly increasing positive integers, and let $M$
be such that relation (\ref{1est}) holds.
In the case $X=R_{net}(P)$ or $X=R_{log}(P)$,
for each $1\leq m\leq M$, let
$$
s_2(m):=\sum_{i=1}^{N_m}\Big [X(i/N_m)- X((i-1)/N_m)\Big ]^2.
$$
In relations (\ref{2est}) and (\ref{3est}), for each $1\leq m\leq M$, 
replacing $N_m/Q(m)$ by $\sqrt{N_m/s_2(m)}$, let
$$
q(m):=\log_2\sqrt{N_m/s_2(m)}\quad\mbox{and}\quad
r(m):=\frac{\log (1/N_m)}{\log \sqrt{s_2(m)/N_m}}
=\frac{\log_2 N_m}{q(m)}.
$$
Then the {\em naive Gladyshev estimator} of the $p$-variation index
is $\gnaive :=\gnaive (X)=r(M)$.
The OLS Gladyshev estimator $\gols$ is defined by
$$
\gols:=\gols (X)=\frac{\sum_{m=1}^M(x_m-\bar x)
\log_2 N_m}{\sum_{m=1}^M(x_m-\bar x)^2},
$$
where $x_m=q(m)=\log_2\sqrt{N_m/s_2(m)}$ and $\bar x=M^{-1}\sum_{m=1}^Mx_m$.

\begin{table}[tb]
\noindent
\caption{\label{summary} OS and G estimation of returns of exchanges rates
based on $\{2^m\colon\,1\leq m\leq 14\}$.} 
\medskip
\centering
\begin{tabular}{|l|c|c|c|c|c|c|c|c|}\hline
Currency &\multicolumn{4}{c|}{net-returns}
         &\multicolumn{4}{c|}{log-returns}\\ \cline{2-9}
    &\rule[0mm]{0mm}{5mm}$\naive$&$\ols$&$\gnaive$&$\gols$&$\naive$&$\ols$&
$\gnaive$&$\gols$ \\  [0.5ex] \hline\hline
AUD/USD & 1.2452 & 1.7617 & 1.3563 & 2.2436 
        & 1.2452 & 1.7664 & 1.3563 & 2.2663\\ \hline
CAD/USD & 1.1493 & 1.8362 & 1.2362 & 2.7382 
        & 1.1493 & 1.8357 & 1.2362 & 2.7422 \\ \hline \hline
DEM/ESP & 1.1504 & 1.7197 & 1.4043 & 2.5684
        & 1.1504 & 1.7274 & 1.4042 & 3.4243 \\ \hline 
DEM/FIM & 1.2124 & 1.8537 & 1.3566 & 2.4970
        & 1.2124 & 1.8571 & 1.3566 & 2.5047 \\ \hline 
DEM/ITL & 1.2420 & 1.7619 & 1.3673 & 2.0673
        & 1.2420 & 1.7552 & 1.3673 & 2.0496 \\ \hline 
DEM/JPY & 1.2671 & 1.8045 & 1.3661 & 2.4714
        & 1.2671 & 1.8034 & 1.3662 & 2.4770 \\ \hline 
DEM/SEK & 1.2483 & 1.7601 & 1.3767 & 2.2612
        & 1.2483 & 1.7558 & 1.3767 & 2.2359 \\ \hline \hline
GBP/DEM & 1.2325 & 1.6459 & 1.3380 & 1.9566
        & 1.2325 & 1.6493 & 1.3380 & 1.9678 \\ \hline 
GBP/USD & 1.2398 & 1.7708 & 1.3334 & 2.1373
        & 1.2398 & 1.7725 & 1.3334 & 2.1466 \\ \hline \hline
USD/BEF & 1.2615 & 1.7083 & 1.4745 & 2.4220
        & 1.2615 & 1.7176 & 1.4745 & 2.5022 \\ \hline 
USD/CHF & 1.2796 & 1.7000 & 1.3902 & 2.0817
        & 1.2796 & 1.7044 & 1.3902 & 2.1001 \\ \hline 
USD/DEM & 1.2432 & 1.7646 & 1.3416 & 2.2111
        & 1.2432 & 1.7673 & 1.3416 & 2.2294 \\ \hline
USD/DKK & 1.2896 & 1.5574 & 1.7142 & 2.2990
        & 1.2894 & 1.6096 & 1.7102 & 2.6473 \\ \hline 
USD/ESP & 1.3574 & 1.9885 & 1.5331 & 2.5367
        & 1.3574 & 2.0155 & 1.5331 & 2.7526 \\ \hline 
USD/FIM & 1.3208 & 1.8583 & 1.4520 & 2.2860
        & 1.3208 & 1.8682 & 1.4520 & 2.3226 \\ \hline 
USD/FRF & 1.2467 & 1.7718 & 1.3643 & 2.2542
        & 1.2467 & 1.7766 & 1.3643 & 2.2811 \\ \hline 
USD/ITL & 1.2977 & 2.1020 & 1.4115 & 3.0934
        & 1.2977 & 2.0863 & 1.4115 & 2.9658 \\ \hline 
USD/NLG & 1.2591 & 1.7570 & 1.3815 & 2.2647
        & 1.2591 & 1.7622 & 1.3816 & 2.2937\\ \hline 
USD/SEK & 1.3263 & 1.9542 & 1.4453 & 2.5304
        & 1.3263 & 1.9606 & 1.4453 & 2.5575 \\ \hline 
USD/XEU & 1.2450 & 1.8503 & 1.3594 & 2.5097
        & 1.2450 & 1.8539 & 1.3594 & 2.5423 \\ \hline \hline
USD/JPY & 1.2703 & 1.8104 & 1.3573 & 2.2434
        & 1.2703 & 1.8143 & 1.3573 & 2.2662 \\ \hline 
USD/MYR & 1.1121 & 1.6715 & 1.3957 & 2.9137
        & 1.1121 & 1.6905 & 1.3956 & 3.0156 \\ \hline 
USD/SGD & 1.1377 & 1.8821 & 1.2584 & 3.2629 
        & 1.1377 & 1.8794 & 1.2584 & 3.2128 \\ \hline 
USD/ZAR & 1.2377 & 1.5460 & 1.4400 & 1.9238
        & 1.2377 & 1.5516 & 1.4401 & 1.9449\\ \hline 
\end{tabular}
\end{table}

%\clearpage

\begin{table}[tb]
\noindent
\caption{\label{metalsum}  OS and G estimation of returns of metal rates
based on $\{2^m\colon\,1\leq m\leq 14\}$.}
\medskip
\centering
\begin{tabular}{|l|c|c|c|c|c|c|c|c|}\hline
Metal   &\multicolumn{4}{c|}{net-returns}
        &\multicolumn{4}{c|}{log-returns}\\ \cline{2-9}
   &\rule[0mm]{0mm}{5mm} $\naive$&$\ols$&$\gnaive$&$\gols$
   &$\naive$&$\ols$&$\gnaive$&$\gols$ \\  [0.5ex] \hline\hline
Gold    & 1.2190 & 1.7249 & 1.3488 & 2.2808 
        & 1.2190 & 1.7216 & 1.3489 & 2.2700 \\ \hline 
Silver  & 1.4097 & 1.8482 & 1.6227 & 2.6066 
        & 1.4097 & 1.8263 & 1.6229 & 2.5728 \\ \hline \hline
Palladium & 1.3371 & 1.6650 & 1.5481 & 2.3535
          & 1.3371 & 1.6603 & 1.5480 & 2.3390 \\ \hline 
Platinum  & 1.2358 & 1.6472 & 1.4134 & 2.2933
          & 1.2358 & 1.6433 & 1.4135 & 2.2790 \\ \hline 
\end{tabular}
\end{table}

\begin{table}[tb]
\noindent
\caption{\label{indexsum}  OS and G estimation of returns of stock indices
based on $\{2^m\colon\,1\leq m\leq 14\}$.}
\medskip
\centering
\begin{tabular}{|l|c|c|c|c|c|c|c|c|}\hline
Index &\multicolumn{4}{c|}{net-returns}
      &\multicolumn{4}{c|}{log-returns}\\ \cline{2-9}
  &\rule[0mm]{0mm}{5mm} $\naive$&$\ols$&$\gnaive$&$\gols$ 
  &$\naive$&$\ols$&$\gnaive$&$\gols$\\  [0.5ex] \hline\hline
SP 500     & 1.1685 & 1.3723 & 1.3726 & 1.9085  
           & 1.1685 & 1.3734 & 1.3727 & 1.9162 \\ \hline
DOW JONES  & 1.1806 & 1.3909 & 1.3937 & 1.9035  
           & 1.1806 & 1.3928 & 1.3938 & 1.9133 \\ \hline 
\end{tabular}
\end{table}

\paragraph*{Estimation results.}
Estimation results for the data set HFDF96 are given by Tables \ref{summary},
\ref{metalsum} and \ref{indexsum}.
More specifically, the tables contain the estimated $p$-variation indices
for the two returns of the bid price associated with the nearest
prior datum.
We picked the log-returns of USD/JPY exchange rates to illustrate 
by Figures \ref{usdjpy1} - \ref{usdjpy4} a difference between the OS and
G estimation results in more detail. 
Figures \ref{usdjpy2} and \ref{usdjpy4} show the estimation results based
on a sequence $\eta=\{N_m=m\colon\,m\geq 1\}$ truncated at $M=240$.
The associated sequence of partitions $\{\lambda (m)\colon\,m\geq 1\}$ in 
this case is not nested, and technical calculations in this case are more 
complex.
Figures \ref{usdjpy1} and \ref{usdjpy3} show the estimation results
based on dyadic partitions, which are used for the rest of results.
Because the conditions of Theorem \ref{main} are more general
as compared to the conditions of the main result of Gladyshev (1961),
the results of the OS estimation are more reliable than the results
of G estimation.
This is also seen from the Figures \ref{usdjpy2} and \ref{usdjpy4}.
The columns $\naive$ and $\ols$ of Tables \ref{summary}, \ref{metalsum} 
and \ref{indexsum} suggest that estimated $p$-variation indices of the 
returns of the financial data are more likely to belong to the interval 
$(1,2)$. 
These columns are essentially the same for net-returns as well as
for log-returns, which would be in agreement if a stock price has the 
quadratic variation along a nested sequence of dyadic partitions 
and its $p$-variation index is not less than 1.
Due to the pattern exhibit by the results of Section \ref{simul:fbm}
when estimating $1/H$, it is unlikely that a fractional Brownian
motion $B_H$ would give a satisfactory fit to the financial data for some 
values of $H$.
However, as seen in Figures \ref{usdjpy1} - \ref{usdjpy4}, a more general 
process from Example \ref{exmp2} might be used to model the HFDF96 data sets.
Especially this concerns a modelling of stock indices 
(see Figure \ref{indexsum}).
Recall that the G estimation is reliable  when applied to sample functions 
of stochastic processes having a suitable relationship between a sample 
function behavior and an asymptotic behavior of the incremental variance.
The estimation results do not reject the hypothesis that returns 
may be modeled by a L\'evy process. 
To be more specific about a degree of data fitting, we need a theoretical
asymptotic analysis of both, the OS and G estimations, which is not
available at this writing.

\paragraph*{Related results and techniques.}
The OS and G estimators are based on properties of a function similar
to a kind of self-similarity property with respect to shrinking 
time intervals, sometimes refered to as a fractal or empirical scaling law.
It is natural that a high-frequency data have already been used to
detect such laws if exist.   
The work of M\"uller et.\ al.\ (1995) addresses this question and provide
some preliminary findings in their analysis of a high-frequency FX data
collected from raw data vendors such as Reuters, Knight-Rider and Telerate. 
Another related work of Mandelbrot (1997) have already been discussed
in Section 4.3 of Norvai\v sa and Salopek (2000).

\section{Discussion}

Parameter estimation of a financial model is a typical econometric 
analysis task.
What is atypical in the preceding analysis is a generality of
the underlying financial model.
This model applies far beyond of limits imposed by the semimartingale theory, 
and its outline can be found in Bick and Willinger (1994), and 
Norvai\v sa (2000a).
The $p$-variation index considered as a parameter of a model is defined
for any function.
Its estimate provides a grade for each concrete continuous time
model of a price process $P$ governed by an exponential, or by a linear 
integral equation having the indefinite integral
\beq\label{prod-int}
P(t)=\lim_{\kappa}\prod_{i=1}^n\big (1+X(t_i\wedge t)-X(t_{i-1}\wedge t)),
\qquad 0<t\leq 1,
\eeq
as its unique solution, where the limit is understood either in the sense of
refinements of partitions $\kappa=\{t_i\colon\,i=0,\dots,n\}$ of $[0,1]$,
or in the more general sense along a fixed sequence of nested partitions.
The linear It\^o stochastic integral equation with respect to
a semimartingale $X$ is one such example of a continuous time model.
A further generality of the underlying financial model could be achieved 
once the net-returns are modified so as to reverse a solution of a 
non-linear integral equation.

\begin{figure}
\onelinecaptionsfalse
\unitlength1cm
\centering
\begin{minipage}[b]{7.0cm}
 \centering
 \includegraphics[height=6.7cm, width=6.7cm]{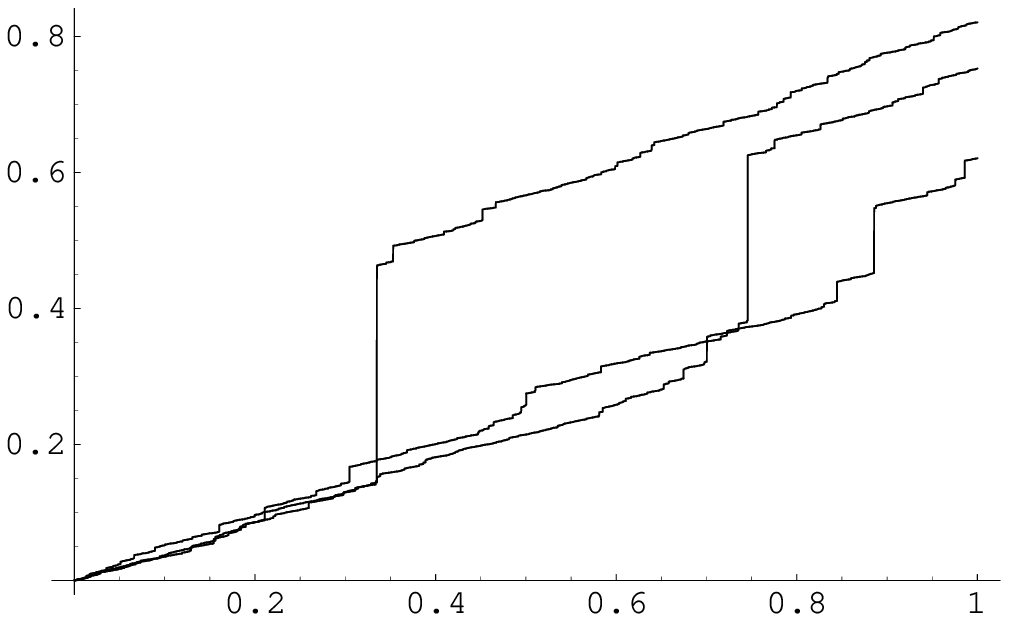}
\end{minipage}%
\hfill
\begin{minipage}[b]{7.0cm}
 \centering
 \includegraphics[height=6.7cm, width=6.7cm]{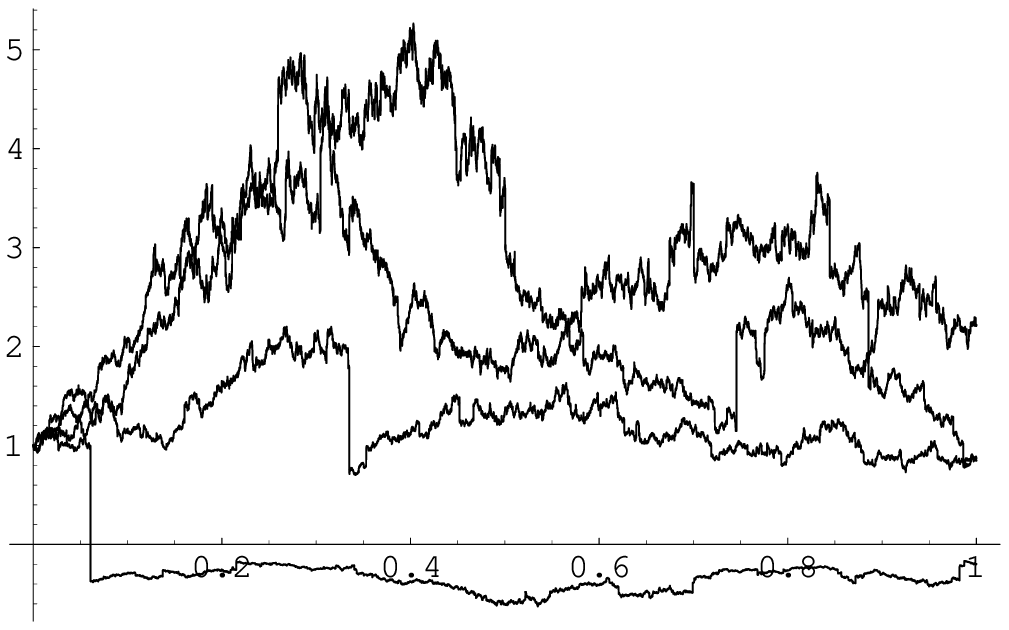}
\end{minipage}\\ 
%\begin{minipage}[t]{7.0cm}
%\captionstyle{flushleft}
%\setcaptionwidth{7.0cm}
\caption{Four trajectories of $R_{net}(P)-R_{log}(P)$ (left) defined by
(\ref{dif_returns}) corresponding to the Dolean exponential 
$P={\cal E}(X_{\alpha})$ (right) given by (\ref{dolean}), where
$X_{\alpha}$ is a S$\alpha$S process with $\alpha=1.7$.
Vertical lines join jumps.}
\label{stable_returns}
%\end{minipage}% 
%\hfill
%\begin{minipage}[t]{7.0cm}
%\setcaptionwidth{7.0cm}
%\caption{OLS estimators: $\ols=1.814$ and $\gols=2.266$
%where $\eta=\{N_m=2^m\colon\,m=1,\dots,14\}$.}
%\label{usdjpy4}
%\end{minipage}%
\end{figure}

We stress the importance of the notion of a return because it
provides a two direction link between theory and practice.  
A financial model without the notion of a return is just an exercise
in theory building.
As we noted earlier, the results of the preceding section
show that the estimated $p$-variation indices for the calculated 
net-returns and log-returns are almost the same.
Here we discuss the difference between the net- and log-returns
$R_{net}(P)-R_{log}(P)$ given by (\ref{dif_returns}) for a simulated 
price process $P$.
Suppose that $P$ is the Dolean exponential ${\cal E}(X_{\alpha})$
of a symmetric $\alpha$-stable process $X_{\alpha}$:
\beq\label{dolean}
{\cal E}(X_{\alpha})(t)
:=\exp\big\{X_{\alpha}(t)-X_{\alpha}(0)\big\}\prod_{(0,t]}
(1+\Delta^{-}X_{\alpha})e^{-\Delta^{-}X_{\alpha}},\qquad
0<t\leq 1.
\eeq   
That is, $P$ is defined by (\ref{prod-int}) with $X$ replaced by $X_{\alpha}$. 
Then the continuous part of the quadratic variation $[P]^c\equiv 0$, and
so the simulation gives the remaining sums in (\ref{dif_returns})
as shown in Figure \ref{stable_returns}.

To simulate the Dolean exponential ${\cal E}(X_{\alpha})$ we use its 
representation as the indefinite product integral (\ref{prod-int}), 
proved by Dudley and Norvai\v sa (1999, Corollary 5.23 in Part II). 
Theoretically, the Dolean exponential composed with the net-return
gives $X_{\alpha}=R_{net}(P)$.
Thus the maximal discrepancy $d:=\sup_t| 
\widetilde{X}_{\alpha}(t)-\widetilde{R}_{net}(\widetilde{P})(t)|$ between 
simulated versions of the two sides should be small if simulation is 
accurate enough.
Indeed, we get that the discrepency $d=3.1\times 10^{-14}$ holds
uniformly for all trajectories in the left Figure \ref{stable_returns}. 
Notice that one of the four trajectories breaks because its Dolean
exponential jumps to the negative side, and so its logarithm is
undefined.

\section*{Appendix A}

This section contains proofs of Theorem \ref{main} 
and Proposition \ref{3main}.
The proofs consist of combining a few relations
between the metric entropy index, the oscillation $\eta$-summing
index and the $p$-variation index, and will be discussed first.

As noted earlier, in the actual calculations of the metric entropy 
it is often simpler to replace closed balls by squares (boxes), 
which leads to the box-counting dimension (see Section 6.1 in 
Cutler, 1993).
First we modify the box-counting method for sets which arise when
connecting a graph of a discontinuous function.
Throughout this section, $f$ is a regulated function defined on a
closed interval $[0,T]$.
This means that for each $t\in (0,T]$, there exists a limit 
$f(t-):=\lim_{u\uparrow t}f(u)$, and for each $t\in [0,T)$, there exists 
a limit $f(t+):=\lim_{u\downarrow t}f(u)$.
For $t\in (0,T]$, let $V_{-}(f(t))$ be the interval connecting the points 
$(t,f(t-))$ and $(t,f(t))$, and for $t\in [0,T)$, let $V_{+}(f(t))$ be the 
interval connecting the points $(t,f(t))$ and $(t,f(t+))$.
On a plane $\RR^2$ with the Cartesian coordinates $(t,x)$, for each
$0\leq u <v\leq T$, let
\beq\label{G_f}
G_f([u,v]):=(u,V_{+}(f(u)))\bigcup\Big (\bigcup\limits_{u< t <v}(t,V_{+}(f(t))
)\bigcup\, (t,V_{-}(f(t)))\Big )\bigcup\, (v,V_{-}(f(v)))\subset\RR^2,
\eeq
where $(t,A(t)):=\{(t,x)\colon\,x\in A(t)\}$.
Let $G_f:=G_f([0,T])$, and let $\gr (f)$ be the graph of the function $f$.
It is clear that $\gr (f)\subset G_f$, and the equality holds between the
two sets if and only if $f$ is continuous on $[0,T]$.
For $\epsilon >0$, a {\em grid} $\CAC_{\epsilon}$ of side length $\epsilon$ 
covering a bounded set $E$, is called a collection of disjoint squares
of equal side length $\epsilon$ whose union is a square containing $E$.
We can also assume that each grid $\CAC_{\epsilon}$ intersects the
origin of the Cartesian coordinates, and its squares have sides parallel
to the coordinates axes.
The usual box-counting method counts all squares of the grid $\CAC_{\epsilon}$
which have nonempty intersection with $E$.
We need a special counting rule when a vertical segment of the set
$E=G_f$ is a part of a vertical line of a grid $\CAC_{\epsilon}$.
Let $\{t_i\colon\,i=0,\dots,N_{\epsilon}\}$ be the partition of $[0,T]$ 
induced by intersecting a grid $\CAC_{\epsilon}$ with the $t$-axes, 
and for each $i=1,\dots,N_{\epsilon}$, let $n_i(\epsilon)$ be the number 
of squares of the grid $\CAC_{\epsilon}$ contained in the strip 
$[t_{i-1},t_i]\times\RR$ which have a nonempty intersection with 
$G_f([t_{i-1},t_i])$.
For each $\epsilon >0$, we then define 
$$
M(G_f;\epsilon):=\sum_{i=1}^{N_{\epsilon}}n_i(\epsilon).
$$
This definition avoids a double covering of certain parts of vertical 
segments, and it reduces to the usual definition if $f$ is continuous.
The present extension allows us to apply Theorem \ref{main} to sample
functions of L\'evy processes, which may have left-side discontinuities.

Next statement shows the usefulness of the box-counting method.

\begin{lem}\label{boxes}
Let $f$ be a regulated function on $[0,T]$. 
For lower and upper metric entropy indices, we have
$$
\Delta^{-}(G_f)=\liminf_{\epsilon\downarrow 0}\frac{\log M(G_f;\epsilon)}
{\log (1/\epsilon)}\quad\mbox{and}\quad
\Delta^{+}(G_f)=\limsup_{\epsilon\downarrow 0}\frac{\log M(G_f;\epsilon)}
{\log (1/\epsilon)}.
$$
\end{lem}

\begin{proof}
For $\epsilon >0$, let $\CAC_{\epsilon}$ be a grid of side length
$\epsilon$ covering $G_f$.
Every square of side $\epsilon$ is included in the ball of diameter
$\epsilon\sqrt{2}$, which, in turn, is included in at most $9$
squares of the grid $\CAC_{\epsilon}$.
Thus
$$
N(G_f;\epsilon\sqrt{2})\leq M(G_f;\epsilon)\leq 9N(G_f;\epsilon\sqrt{2}),
$$
proving the lemma.
\qed\end{proof}

Next we look at replacing a limit as $\epsilon\downarrow 0$ by a limit
along a countable sequence.
Let $\{\epsilon_m\colon\,m\geq 1\}$ be a sequence of positive real 
numbers strictly decreasing to $0$, and for each $m\geq 1$, 
let $\Delta_m:=(\epsilon_{m+1},\epsilon_m]$.
For a family $\{A_{\epsilon}\colon\,\epsilon >0\}$ of numbers,  we have
$$
\limsup_{\epsilon\downarrow 0}A_{\epsilon}=\limsup_{m\to\infty}
\sup_{\epsilon\in\Delta_m}A_{\epsilon}\geq\limsup_{m\to\infty}
A_{\epsilon_m}
$$
and
$$
\liminf_{\epsilon\downarrow 0}A_{\epsilon}=\liminf_{m\to\infty}
\inf_{\epsilon\in\Delta_m}A_{\epsilon}\leq\liminf_{m\to\infty}
A_{\epsilon_m}.
$$
Thus by the preceding lemma, we have that
\beq\label{2subsequences}
\Delta^{-}(G_f)\leq\liminf_{m\to\infty}\frac{\log M(G_f;\epsilon_m)}
{\log (1/\epsilon_m)}\quad\mbox{and}\quad
\limsup_{m\to\infty}\frac{\log M(G_f;\epsilon_m)}{\log (1/\epsilon_m)}
\leq \Delta^{+}(G_f).
\eeq

\begin{lem}\label{subsequences}
Let $f$ be a regulated function on $[0,T]$, and let $\{\epsilon_m
\colon\,m\geq 1\}$ be a sequence of positive real numbers strictly
decreasing to $0$ so that
\beq\label{1subsequences}
\lim_{m\to\infty}\frac{\log \epsilon_m}{\log\epsilon_{m+1}}=1.
\eeq
Then
\beq\label{3subsequences}
\left\{ \begin{array}{l}
\liminf_{\epsilon\downarrow 0}\\
\limsup_{\epsilon\uparrow 0}
\end{array} \right\}
\frac{\log M (G_f;\epsilon)}{\log (1/\epsilon)}
=\left\{ \begin{array}{l}
\liminf_{m\to\infty}\\
\limsup_{m\to\infty}
\end{array} \right\}
\frac{\log M (G_f;\epsilon_m)}{\log (1/\epsilon_m)}.
\eeq
\end{lem}

\begin{proof}
It is enough to prove the reverse inequalities in relations
(\ref{2subsequences}). 
For $0<\epsilon\leq\epsilon_1$, let $m\geq 1$ be such that 
$\epsilon\in\Delta_m=(\epsilon_{m+1},\epsilon_m]$.
Since a square with a side length $\epsilon$ is contained in at most
four squares with a side length $\epsilon_m$, we have
$ M(G_f;\epsilon_m)\leq 4 M (G_f;\epsilon)$.
Similarly, $M(G_f;\epsilon)\leq 4 M (G_f;\epsilon_{m+1})$.
Thus for each $\epsilon\in\Delta_m$,
$$
\frac{\log \epsilon_m}{\log\epsilon_{m+1}}\Big (\frac{\log M(G_f;\epsilon_m)
-\log 4}{\log 1/\epsilon_m}\Big )\leq 
\frac{\log M(G_f;\epsilon)}{\log (1/\epsilon )}
\leq\frac{\log \epsilon_{m+1}}{\log\epsilon_{m}}\Big (\frac{\log M
(G_f;\epsilon_{m+1})+\log 4}{\log 1/\epsilon_{m+1}}\Big ).
$$
The conclusion now follows by the assumption (\ref{1subsequences}).
\qed\end{proof}

Now recall relation (\ref{Q(f)}) defining the oscillation $\eta$-summing 
sequence $\{Q(f;\lambda (m))\colon\,m\geq 1\}$.

\begin{lem}\label{CCM}
Let $f$ be a regulated non-constant function on $[0,T]$.
For any sequence $\{N_m\colon\,m\geq 1\}$ of strictly increasing positive 
integers, we have
\beq\label{1CCM}
\left\{\begin{array}{l}
\liminf_{m\to\infty} \\ \limsup_{m\to\infty}
\end{array}\right\}
\frac{\log M (G_f;T/N_m)}{\log N_m} 
=\left\{\begin{array}{l}
\liminf_{m\to\infty} \\ \limsup_{m\to\infty}
\end{array}\right\}
\frac{\log N_m Q (f;\lambda (m))}{\log N_m}.
\eeq
\end{lem}

\begin{proof}
Let $\{N_m\colon\,m\geq 1\}$ be a sequence of strictly increasing positive 
integers.
For each $m\geq 1$, let $A_m$ and $B_m$ be the terms under the limit
signs in equality (\ref{1CCM}) in the given order.
It is enough to prove that
\beq\label{3CCM}
\liminf_mA_m\leq \liminf_mB_m,\qquad
\limsup_m A_m\leq\limsup_mB_m
\eeq
and 
\beq\label{4CCM}
\liminf_mB_m\leq \liminf_mA_m,\qquad
\limsup_m B_m\leq\limsup_mA_m.
\eeq
Suppose that the set $G_f$ is on a plane with the Cartesian coordinates 
$(t,x)$.
For each integer $m\geq 1$, let $\CAC_m$ be the grid of side length $T/N_m$ 
intersecting the origin and with sides parallel to the coordinate axes.
Thus $\CAC_m$ intersects with the $t$-axes at each point of the partition
$\lambda (m)=\{t_i^m=iT/N_m\colon\,i=0,\dots,N_m\}$.
Let $m\geq 1$.
For $i=1,\dots,N_m$, let $n_i=n_i(m)$ be the number of squares of the grid
$\CAC_m$ contained in the strip $\Delta_{i,m}\times\RR$ and
covering the set $G_f(\Delta_{i,m})$ defined by relation (\ref{G_f}),
where $\Delta_{i,m}=[t_{i-1}^m,t_i^m]$.
By the definition of $G_f(\Delta_{i,m})$, 
$$
(N_m/T)\osc (f;\Delta_{i,m})\leq n_i\leq 2+(N_m/T)\osc (f;\Delta_{i,m}),
$$ 
for each $i=1,\dots,N_m$.
Summing over all indices $i$, it follows that the inequalities
\beq\label{2CCM}
M(G_f;T/N_m)-2N_m\leq (N_m/T)Q(f;\lambda (m))\leq M(G_f;T/N_m)
\eeq
hold for each $m\geq 1$.
By the second inequality in display (\ref{2CCM}), we have that 
inequalities (\ref{4CCM}). 
To prove inequalities (\ref{3CCM}), first suppose that $\liminf_mB_m=1$.
Then by relation (\ref{Qin[1,2]}), $\limsup_mB_m=1$, and inequalities 
(\ref{3CCM}) follow because $\liminf_mA_m\geq 1$.
Now suppose that $\liminf_mB_m>1$.
In that case, inequalities (\ref{3CCM}) follow from
the first inequality in display (\ref{2CCM}) and from the following two 
relations:
$$
\liminf_{m\to\infty}\frac{\log a_m}{\log b_m}
=\sup\big\{\alpha >0\colon\,\lim_{m\to\infty}b_m^{-\alpha}a_m=+\infty\big\}
=\sup\big\{\alpha >0\colon\,\inf_{m\geq 1}b_m^{-\alpha}a_m>0\big\}
$$
and
\beq\label{s-limsup}
\limsup_{m\to\infty}\frac{\log a_m}{\log b_m}
=\inf\big\{\alpha >0\colon\,\lim_{m\to\infty}b_m^{-\alpha}a_m=0\big\}
=\inf\big\{\alpha >0\colon\,\sup_{m\geq 1}b_m^{-\alpha}a_m<+\infty\big\},
\eeq
valid for any two sequences $\{a_m\colon\,m\geq 1\}$ and $\{b_m\colon\,
m\geq 1\}$ of positive numbers such that $\lim_{m\to\infty}b_m=+\infty$.
\qed\end{proof}

Combining Lemmas \ref{boxes}, \ref{subsequences} and \ref{CCM}, 
it follows that the following statement holds.

\begin{cor}
Let $f$ be a regulated non-constant function on $[0,T]$, and let 
$\{N_m\colon\,m\geq 1\}$ be a sequence of strictly increasing positive 
integers such that
\beq\label{4subsequences}
\lim_{m\to\infty}\frac{\log N_m}{\log N_{m+1}}=1.
\eeq
Then
$$
\Delta (G_f)=\lim_{\epsilon\downarrow 0}\frac{\log M(G_f;\epsilon)}{\log
(1/\epsilon )}=\lim_{m\to\infty}\frac{\log N_mQ(f;\lambda (m))}{\log N_m},
$$
provided that at least one of the three limits exists and is finite.
\end{cor}

Next is the final step in a chain of inequalities used to
prove Theorem \ref{main}.

\begin{lem}\label{bc-pv}
Let $f$ be a regulated function on $[0,T]$, and let 
$\eta=\{N_m\colon\,m\geq 1\}$ be a sequence of strictly increasing 
positive integers.
Then
\beq\label{1bc-pv}
\Big (2-\delta_{\eta}^{+}(f)=\Big )\qquad
\limsup_{m\to\infty}\frac{\log N_mQ(f;\lambda (m))}{\log N_m}\leq 
2-\frac{1}{1\vee\upsilon (f)}.
\eeq
\end{lem}

\begin{proof}
If $\upsilon (f)<1$ then $f$ has bounded $1$-variation, and so
(\ref{1bc-pv}) holds by relation (\ref{bv}).
If $\upsilon (f)=+\infty$, then (\ref{1bc-pv}) holds because
$\delta_{\eta}^{+}(f)\geq 0$ by relation (\ref{Qin[1,2]}).
Thus one can assume that $1\leq \upsilon (f)<+\infty$.
Let $v_p(f;[0,T])<\infty$ for some $1\leq p<\infty$.
We claim that for each integer $m\geq 1$,
\beq\label{3bc-pv}
N_m^{\frac{1}{p}-1}Q(f;\lambda (m))\leq v_p(f;[0,T])^{1/p}.
\eeq
Indeed, if $p=1$ then $Q(f;\lambda (m))\leq v_1(f;[0,T])$.
Suppose that $p>1$ and $m\geq 1$.
Since for each $i=1,\dots,N_m$, $\osc (f;\Delta_{i,m})\leq v_p(f;
\Delta_{i,m})^{1/p}$, by H\"older's inequality, we have
\begin{eqnarray*}
Q (f;\lambda (m))&\leq& \sum_{i=1}^{N_m}v_p(f;\Delta_{i,m})^{1/p}
\leq N_m^{1-1/p}\Big (\sum_{i=1}^{N_m}v_p(f;\Delta_{i,m})\Big )^{1/p}\\[2mm]
\mbox{by subadditivity of $v_p$}\quad &\leq & N_m^{1-1/p}v_p(f;[0,T])^{1/p},
\end{eqnarray*}
proving relation (\ref{3bc-pv}).
If $p>\upsilon (f)\geq 1$, then by (\ref{3bc-pv}), 
$$
\sup_{m\geq 1}N_m^{-(2-1/p)}\big [N_mQ(f;\lambda (m))\big ]
\leq v_p(f;[0,T])^{1/p}<+\infty.
$$
Thus relation (\ref{1bc-pv}) follows by the relation (\ref{s-limsup}), 
proving the lemma.
\qed\end{proof}

Now we are ready to complete the proofs.

\medskip
\noindent
{\bf Proof of Theorem \ref{main}.}
Let $\eta=\{N_m\colon\,m\geq 1\}$ be a sequence of strictly increasing
positive integers.
We have to prove that $\delta_{\eta}(f)=1/(1\vee\upsilon (f))$.
By Definition \ref{OS-index}, this will be done once we will
show that
\beq\label{1main}
\lim_{m\to\infty}\frac{\log N_m Q(f;\lambda (m))}{\log N_m}
=2-\frac{1}{1\vee\upsilon (f)}. 
\eeq
Since the graph $\gr (f)\subset G_f$, we have
\begin{eqnarray*}
\Delta^{-}(\gr (f)) &\leq & \Delta^{-}(G_f)=\liminf_{
\epsilon\downarrow 0}\frac{\log M(G_f;\epsilon)}{\log (1/\epsilon)}
\quad\stackrel{\mbox{by (\ref{2subsequences})}}{\leq}\quad
\liminf_{m\to\infty}\frac{\log M(G_f;1/N_m)}{\log N_m}\\
\mbox{by (\ref{1CCM})}\quad
&=&\liminf_{m\to\infty}\frac{\log N_mQ(f;\lambda (m))}{\log N_m}
\leq \limsup_{m\to\infty}\frac{\log N_mQ(f;\lambda (m))}{\log N_m}\\
\mbox{by (\ref{1bc-pv})}\quad
&\leq& 2-\frac{1}{1\vee\upsilon (f)}.
\end{eqnarray*}
Since the left and right sides are equal by assumption (\ref{2main}),
relation (\ref{1main}) holds, proving Theorem \ref{main}.
\qed

\medskip
\noindent
{\bf Proof of Proposition \ref{3main}.}
If $f$ is constant then the conclusion clearly holds,
and so we can assume that $f$ is non-constant.
Let $\eta=\{N_m\colon\,m\geq 1\}$ be a sequence of strictly increasing
positive integers such that relation (\ref{4subsequences}) holds.
Since the graph $\gr (f)\subset G_f$, we have
\begin{eqnarray*}
\Delta^{+}(\gr (f)) & =&\limsup_{\epsilon\downarrow 0}
\frac{\log M(G_f;\epsilon)}{\log (1/\epsilon)}
\quad\stackrel{\mbox{by (\ref{3subsequences})}}{=}\quad
\limsup_{m\to\infty}\frac{\log M(G_f;1/N_m)}{\log N_m}\\
\mbox{by (\ref{1CCM})}\quad
&=&\limsup_{m\to\infty}\frac{\log N_mQ(f;\lambda (m))}{\log N_m}
\quad\stackrel{\mbox{by (\ref{1bc-pv})}}{\leq}\quad
2-\frac{1}{1\vee\upsilon (f)},
\end{eqnarray*}
proving Proposition \ref{3main}.
\qed

%\end{appendix}

\vspace*{0.1truein}
%\addcontentsline{toc}{chapter}{\protect\numberline{}{Bibliography}}

\section*{References}

\parindent=2pc\def\hang#1{\vbox{\hsize\the\hsize\noindent%
     \hangindent\the\parindent#1\par}\smallskip}%

\hang{S.\ M.\ Berman,
``Harmonic analysis of local times and sample functions of
Gaussian processes,''
{\em Trans.\ Amer.\ Math.\ Soc.} vol.\ 143 pp.\ 269-281, 1969.}

\hang{A.\ Bick and W.\ Willinger,
``Dynamic spanning without probabilities,''
{\em Stoch.\ Proc.\ Appl.} vol.\ 50 pp.\ 349-374, 1994.} 

\hang{R.\ M.\ Blumenthal and R.\ K.\ Getoor,
``Sample functions of stochastic processes with stationary independent
increments,''
{\em J.\ Math.\ and Mech.} vol.\ 10 pp.\ 493-516, 1961.}

\hang{R.\ M.\ Blumenthal and R.\ K.\ Getoor,
``The dimension of the set of zeros and the graph of a symmetric
stable process,'' 
{\em Illinois J.\ Math.} vol.\ 6 pp.\ 308-316, 1962.}

\hang{J.\ Y.\ Campbell, A.\ W.\ Lo and A.\ C.\ MacKinlay,
{\em The econometrics of financial markets},
Princeton University Press: New Jersey, 1997.}

\hang{P.\ H.\ Carter, R.\ Cawley and R.\ D.\ Mauldin,
``Mathematics of dimension measurement for graphs of functions,''
in {\em Fractal Aspects of Materials: Disordered Systems},
Eds. D.\ A.\ Weitz, L.\ Sander and B.\ Mandelbrot,
Material Research Society, US, 1988, pp.\ 183-186. }

\hang{J.\ M.\ Chambers, C.\ L.\ Mallows and B.\ W.\ Stuck,
``A method for simulating stable random variables,''
{\em J.\ Amer.\ Statist.\ Assoc.} vol.\ 71 no. 354 pp.\ 340-344, 1976.}

\hang{H.\ Cramer and M.\ R.\ Leadbetter,
{\em Stationary and related stochastic processes},
Wiley: New York, 1967.}

\hang{M.\ E.\ Crovella and M.\ S.\ Taqqu,
``Estimating the heavy tail index from scaling properties,''
{\em Meth.\ and Comp.\ in Appl.\ Probab.} vol.\ 1 pp.\ 55-79, 1999.}

\hang{C.\ D.\ Cutler,
``A review of the theory and estimation of fractal dimension,''
in {\em Dimension Estimation and Models}, Ed.\ H.\ Tong.,
World Scientific, Singapore, 1993, pp.\ 1-107.}

\hang{B.\ Dubuc, J.\ F.\ Quiniou, C.\ Roques-Carmes, C.\ Tricot
and S.\ W.\ Zucker,
``Evaluating the fractal dimension of profiles,''
{\em Physical Review A} vol.\ 39 pp.\ 1500-1512, 1989.}

\hang{R.\ M.\ Dudley and R.\ Norvai\v sa,
{\em Differentiability of Six Operators on Nonsmooth Functions and
p-Variation},
Lect.\ Notes in Math., vol.\  1703,
Springer: Berlin, 1999}

\hang{H.\ F\"olmer,
``Calcul d'It\^o sans probabilit\'es,''
in {\em S\'eminaire de Probabilit\'es} XV, Eds.\ J.\ Az\'ema and 
M.\ Yor, Lect.\ Notes Math.\ vol.\ 850, Springer, 1981, pp.\ 143-150.}

\hang{E.\ G.\ Gladyshev,
``A new limit theorem for stochastic processes with Gaussian increments,''
{\em Theor.\ Probability Appl.} vol.\ 6 pp.\ 52-61, 1961.}

\hang{P.\ Hall and A.\ Wood,
``On the performance of box-counting estimators of fractal dimension,''
{\em Biometrika}, vol.\ 80 pp.\ 256-252, 1993.}

\hang{N.\ C.\ Jain and D.\ Monrad,
``Gaussian measures in $B_p$,''
{\em Ann.\ Probab.} vol.\ 11 pp.\ 46-57, 1983.}

\hang{N.\ K\^ono,
``Hausdorff dimension of sample paths for self-similar processes,''
in {\em Dependence in Probability and Statistics}, Eds.,
E.\ Eberlein and M.\ S.\ Taqqu,
Birkh\"auser, Boston, 1986, pp.\ 109-117.}

\hang{R.\ E.\ Maeder,
``Fractional Brownian motion,''
{\em Mathematica J} vol.\ 6 pp.\ 38-48, 1995.}

\hang{B.\ B.\ Mandelbrot,
{\em Fractals and Scaling in Finance: Discontinuity, Concentration,
Risk}, Selecta volume E,
Springer: New York, 1997.}

\hang{I.\ Monroe,
``On the $\gamma$-variation of processes with stationary independent
increments,''
{\em Ann.\ Math.\ Statist.} vol.\ 43 pp.\ 1213-1220, 1972.}

\hang{U.\ A.\ M\"uller, M.\ M.\ Dacorogna, R.\ D.\ Dav\'e, O.\ V.\ Pictet,
R.\ B.\ Olsen and J.\ R.\ Ward,
``Fractals and intrinsic time - a chalenge to econometricians,''
Internal document UAM.1993-08-16, Olsen \& Associates, 
Switzerland, 1995. http://www.olsen.ch}

\hang{R.\ Norvai\v sa,
``Modelling of stock price changes: A real analysis approach,''
{\em Finance and Stochastics} vol.\ 4 pp.\ 343-369, 2000a.}

\hang{R.\ Norvai\v sa,
``Quadratic variation, $p$-variation and integration with applications
to stock price modelling,'' 2000b. (in preparation)}

\hang{R.\ Norvai\v sa and D.\ M.\ Salopek,
``Estimating the Orey index of a Gaussian stochastic process with
stationary increments: An application to financial data set,''
in {\em Stochastic Models}, Proc.\ Int.\ Conf., Ottawa, Canada
June 10-13, 1998, Eds.\ L.\ G.\ Gorostiza and B.\ G.\ Ivanoff.
Canadian Math.\ Soc., Conference Proceedings, vol.\ 26 2000, pp.\ 353-374.}

\hang{R.\ Norvai\v sa and D.\ M.\ Salopek,
``Supplement to the present paper, containing codes in Mathematica,''
2000b. (available from the authors upon request)}

\begin{sloppypar}
\hang{S.\ Orey,
``Gaussian sample functions and the Hausdorff dimension of level crossings,''
{\em Z.\ Wahrsch.\ verw.\ Gebiete} vol.\ 15 pp.\ 249-256, 1970.}
\end{sloppypar}

\hang{W.\ E.\ Pruit and S.\ J.\ Taylor,
``Sample path properties of processes with stable components,''
{\em Z.\ Wahrsch.\ verw.\ Geb.} vol.\ 12 pp.\ 267-289, 1969.}

%\hang{G.\ Samorodnitsky and M.\ Taqqu,}

%\hang{C.\ Tricot,
%{\em Curves and Fractal Dimension}.
%Springer, New York, 1995.}

\hang{V.\ M.\ Zolotarev,
{\em One-dimensional Stable Distributions},
Translations of Mathematical Monographs, vol.\ 65,
American Mathematical Society, 1986.}

\end{document}